# PROBABILISTIC VALIDATION OF HOMOLOGY COMPUTATIONS FOR NODAL DOMAINS

By Konstantin Mischaikow[1] and Thomas Wanner[2]

*Rutgers University and George Mason University*

Homology has long been accepted as an important computable tool for quantifying complex structures. In many applications, these structures arise as nodal domains of real-valued functions and are therefore amenable only to a numerical study based on suitable discretizations. Such an approach immediately raises the question of how accurate the resulting homology computations are. In this paper, we present a probabilistic approach to quantifying the validity of homology computations for nodal domains of random fields in one and two space dimensions, which furnishes explicit probabilistic a priori bounds for the suitability of certain discretization sizes. We illustrate our results for the special cases of random periodic fields and random trigonometric polynomials.

**1. Introduction.** The practical need to extract low-dimensional nonlinear structures from high-dimensional data sets has led to the introduction of topological methods in statistical analysis. There is a growing body of literature [10, 13, 14, 26, 28] that addresses the problem of providing efficient algorithms for estimating topological properties of a fixed, but unknown, manifold $X$ given a point-cloud data set that lies on or near $X$. The equally important complementary question is that of the accuracy of these estimations. To the best of our knowledge, the first result is due to Niyogi, Smale and Weinberger [25]. In this paper, the authors propose a stochastic algorithm for computing the homology of a given manifold $X \subset \mathbb{R}^d$ by randomly sampling $M$ points from the manifold and derive explicit bounds on the probability that their algorithm computes the correct homology. The probability bound depends on the number $M$ and a condition number $1/\tau$. The

Received September 2006; revised December 2006.

[1]Supported in part by NSF Grants DMS-05-11115 and DMS-01-07396, by DARPA and the U.S. Department of Energy.

[2]Supported in part by NSF Grant DMS-04-06231 and the U.S. Department of Energy under Contract DE-FG02-05ER25712.

*[AMS 2000 subject classifications.](http://www.ams.org/msc/)* 60G60, 55N99, 60G15, 60G17.

latter parameter encodes both local curvature information of the manifold $X$ and global separation properties. More precisely, the inverse condition number $\tau$ is the largest number such that the open normal bundle about $X \subset \mathbb{R}^d$ of radius $r$ is embedded in $\mathbb{R}^d$ for all $r < \tau$.

Although stimulated by the aforementioned results, the motivation for the work presented in this paper comes from the study of patterns produced by nonlinear evolutionary systems such as phase separation in materials, fluid flow and population biology. While the mathematical models for such systems are typically infinite-dimensional, understanding the time evolution of the multitude of often complicated patterns produced on an underlying two- or three-dimensional domain is of considerable, if not primary, interest. As is described in greater detail below, algebraic topology and, in particular, computational homology can be used to obtain new statistical measurements for either numerical or experimental data associated with these spatially and temporally complex systems [18, 19, 23]. As above, the accuracy of these homology computations is an important question. However, because the data is obtained through numerical methods or by standard imaging techniques, the data points are not naturally modeled via a random sampling, but rather appear on a fixed grid. On the other hand, the patterns are constantly changing and thus one is interested in the topology of a multitude of manifolds. Furthermore, at those points in time for which the topology of the pattern changes, the condition number $1/\tau$ becomes unbounded.

These observations suggest the need for a complementary result concerning the accuracy of homology computations. To be more precise, consider a compact rectangular domain $G := [a,b]^d \subset \mathbb{R}^d$ and a smooth function $u \in C^2(G,\mathbb{R})$. The *nodal domains* of $u$ are given by the sub- and super-level sets

$$N^\pm := \{x \in G : \pm u(x) \geq 0\},$$

which can be interpreted as patterns on the domain $G$. These are the objects which we wish to identify through their homology, that is, we are interested in the statistics of $H_*(N^\pm)$ as the functions $u$ vary.

Given the limitations of numerical and imaging (think digital camera) techniques, the nodal domains must be approximated. The following provides a reasonable model for such an approximation.

DEFINITION 1.1 (*Cubical approximation of a nodal domain*). Consider a compact rectangular domain $G := [a,b]^d \subset \mathbb{R}^d$ and a continuous function $u : G \to \mathbb{R}$. Let $M$ be an arbitrary positive integer. Define the *equidistant M-discretization* of $[a,b]$ as the collection of grid points

$$x_k := a + k \cdot \frac{b-a}{M}, \qquad \text{for } k = 0, \ldots, M.$$

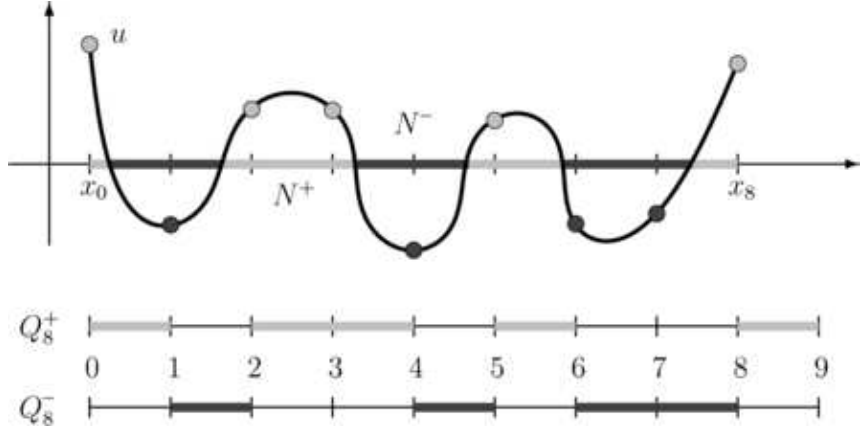

FIG. 1. *Nodal domain approximation using cubical sets in one dimension.*

The *cubical approximations* $Q_M^\pm$ of the nodal domains $N^\pm = \{x \in G : \pm u(x) \geq 0\}$ are defined as the sets

$$Q_M^\pm := \bigcup \left\{ \prod_{\ell=1}^d [k_\ell, k_\ell + 1] : \pm u(x_{1,k_1}, \ldots, x_{d,k_d}) \geq 0, k \in \{0, \ldots, M\}^d \right\},$$

where $k = (k_1, \ldots, k_d)$ and where $x_{j,0}, \ldots, x_{j,M}$ denotes the equidistant $M$-discretization of the $j$th component interval of $G$. For the case of one-dimensional domains, this definition is illustrated in Figure 1; for the two-dimensional case, see Figure 2.

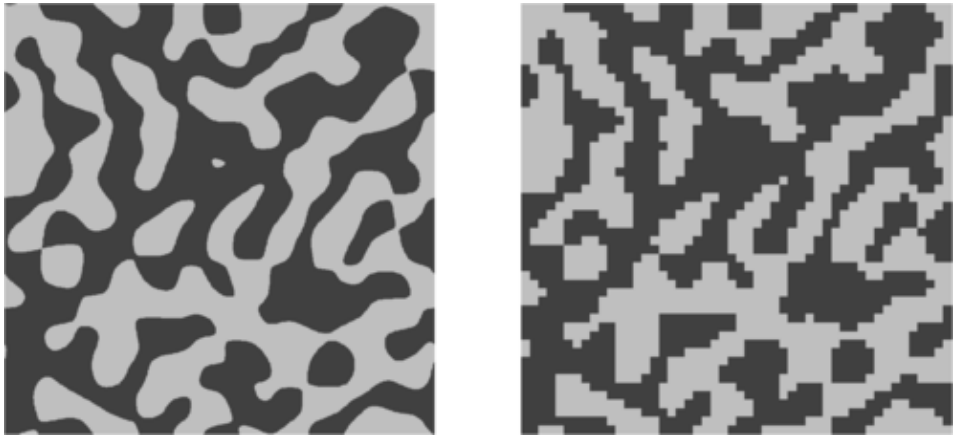

FIG. 2. *Nodal domains of a random trigonometric polynomial in two space dimensions and their cubical approximations with $M = 50$.*

At first glance, one might expect that for sufficiently smooth functions $u$ and sufficiently large discretization size $M$, the homologies of the nodal domains $N^\pm$ and of their cubical approximations $Q_M^\pm$ would coincide. Simple counterexamples to this exist; however, this assertion can be made correct from a probabilistic point of view.

Consider a probability space $(\Omega, \mathcal{F}, \mathbb{P})$ and let $u: G \times \Omega \to \mathbb{R}$ denote a random field over $(\Omega, \mathcal{F}, \mathbb{P})$ such that for $\mathbb{P}$-almost all $\omega \in \Omega$, the function $u(\cdot, \omega): G \to \mathbb{R}$ is twice continuously differentiable. We will impose the following two assumptions throughout this paper:

(A1) for every $x \in G$, we have $\mathbb{P}\{u(x) = 0\} = 0$;
(A2) the random field is such that $\mathbb{P}\{0 \text{ is a critical value of } u\} = 0$.

To simplify the notation, we follow the customary procedure of dropping the explicit $\omega$-dependence in the formulation of events, that is, we write, for example, $\mathbb{P}\{u(x) = 0\}$ instead of $\mathbb{P}(\{\omega \in \Omega : u(x, \omega) = 0\})$.

THEOREM 1.2. *Consider a probability space $(\Omega, \mathcal{F}, \mathbb{P})$, a domain $G := [a,b]^d \subset \mathbb{R}^d$ and a random field $u: G \times \Omega \to \mathbb{R}$ over $(\Omega, \mathcal{F}, \mathbb{P})$ such that for $\mathbb{P}$-almost all $\omega \in \Omega$, the function $u(\cdot, \omega): G \to \mathbb{R}$ is twice continuously differentiable. For each $\omega \in \Omega$, denote the nodal domains of $u(\cdot, \omega)$ by $N^\pm(\omega) \subset G$ and given a positive integer $M$, denote their cubical approximations by $Q_M^\pm(\omega)$.*

*If both (A1) and (A2) are satisfied, then for $\mathbb{P}$-almost all $\omega \in \Omega$, the following holds. For all sufficiently large values of $M$, the homology of the cubical approximations $Q_M^\pm(\omega)$ matches that of the nodal domains $N^\pm(\omega)$. In other words, there exists a random variable $\bar{M}: \Omega \to \mathbb{N}$ such that*

$$\mathbb{P}\{\text{For all } M \geq \bar{M} \text{ one has } H_*(N^\pm) \cong H_*(Q_M^\pm)\} = 1.$$

*Note that the random variable $\bar{M}$ is in general neither constant nor bounded.*

REMARK 1.3. One can show that in the situation of Theorem 1.2, there actually exists a homeomorphism $h: G \to G$ which is isotopic to the identity such that $h(N^\pm) = Q_M^\pm$. In other words, the nodal sets and their cubical approximations are homeomorphic. Similar extensions are also valid for our probabilistic main results. Nevertheless, in view of the current state of applications, we present all of our results within the homology setting.

We do not provide an explicit proof of this result since it follows from the realization that for $\mathbb{P}$-almost all $\omega \in \Omega$, the level set $u(x, \omega) = 0$ is a manifold with finite condition number $1/\tau < \infty$ and a straightforward modification of the proof of [25], Proposition 3.1, then readily furnishes the above result.

Observe that the hypotheses for Theorem 1.2 provide a fairly general framework in which to study the topology of random manifolds. However,

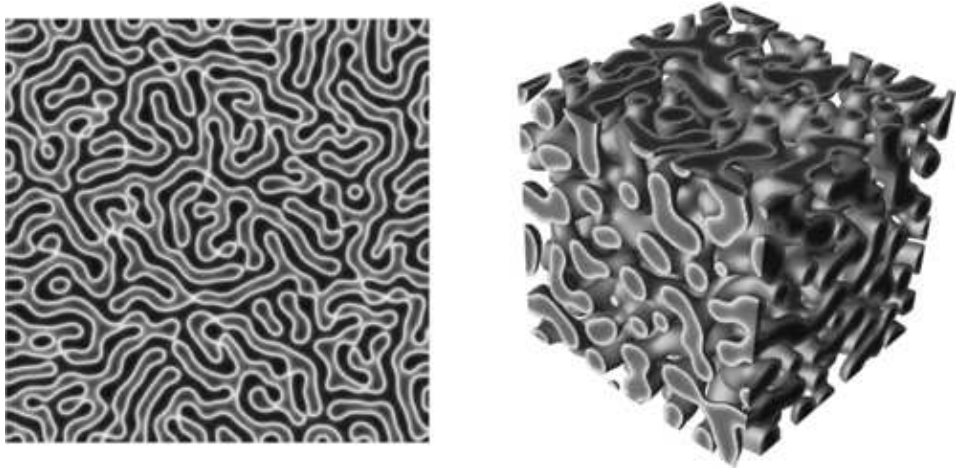

Fig. 3. *Complicated microstructures generated from simulations of the Cahn–Hilliard—Cook model in two and three space dimensions.*

the motivation for this work leads us in a slightly different direction. As an example, consider the problem of phase separation in binary metal alloys. During the last few decades, many models have been proposed for describing this process. Starting from the celebrated Cahn–Hilliard partial differential equation model [8, 9], numerous extensions and refinements have been proposed over the years. While all of these models are infinite-dimensional evolution equations they can be deterministic or stochastic partial differential equations and can contain nonlocal or other additional terms. Moreover, all of these models are phenomenological models in the sense that while they are based on important physical insights, they are not rigorously derived from first principles and so assessing the quality of a specific model becomes a central and crucial problem. Since the material microstructures which are generated by the phase separation process have direct implications for material properties, any such assessment should consider to what extent the complicated patterns by the models match experimental observations.

In the above-described setting of evolution equations, patterns are usually described in a very specific framework. Returning to the materials example, the various models of Cahn–Hilliard type describe the temporal evolution of a phase variable $u = u(t, x)$, where $x \in G$ describes the underlying container for the alloy. The scalar phase variable $u$ describes the local alloy concentration in the following way. If $u$ is close to $+1$, the local concentration consists only of the first alloy component, values close to $-1$ represent the second component and values in between indicate corresponding mixtures of both. Typical patterns that are produced by the Cahn–Hilliard models are shown in Figure 3 and mathematical results have recently been obtained which provide some insight into their dynamics [6, 7, 30]. Since we are interested in regions where one of the two components dominates, the microstructures are generally thought of as nodal domains of the function $u$,

$$N^\pm(t) = \{x \in G : \pm u(t,x) \geq 0\}.$$

How well do these nodal domains match actual experimental results? In a series of papers [20, 21, 24], Hyde et al. described experimentally obtained microstructures in iron-chromium alloys using a variety of techniques. One of their measurements was the handle number of the structure which, from a mathematical point of view corresponds to studying the first Betti number of the nodal domains $N^\pm(t)$ as a function of $t$. In [19], we used methods from computational homology to numerically study the Betti numbers of the nodal domains, in an attempt to unveil differences between two different models, the deterministic Cahn–Hilliard equation and the stochastic Cahn–Hilliard–Cook model. Since the underlying nonlinear evolutions were solved numerically, we do not have direct access to the nodal domains. However, for a fixed grid size bounded from below by the resolution used in the numerical approximation, the associated cubical approximation is easily obtained.

This last comment indicates the limitations of Theorem 1.2. Ideally, given an evolution equation and a fixed cubical approximation of the nodal domain, we would like to be able to estimate the amount of time for which the homology computations of the nodal domains are correct. This problem appears to be extremely difficult. Thus we pose the following, simpler, goal. Consider a random generalized Fourier series of the form

$$(1) \qquad u(x,\omega) = \sum_{k=0}^{\infty} \alpha_k \cdot g_k(\omega) \cdot \varphi_k(x), \qquad u : G \times \Omega \to \mathbb{R},$$

where $G \subset \mathbb{R}^d$ denotes a compact rectangular set and $\varphi_k$, $k \in \mathbb{N}_0$, denotes a sequence of basis functions which will be specified in more detail later. Assume that the $g_k$ are independent identically distributed real-valued random variables over a common probability space $(\Omega, \mathcal{F}, \mathbb{P})$. The coefficients $\alpha_k$ denote deterministic real numbers which will be related to the smoothness properties of the realizations $u(\cdot,\omega) : G \to \mathbb{R}$. Let $M$ be an arbitrary positive integer. Approximate the random nodal domains $N^\pm(\omega)$ of the realization $u(\cdot,\omega)$ by random cubical sets $Q_M^\pm(\omega)$, as in Definition 1.1. Estimating the accuracy of our homology computations then reduces to the following problem:

(P) *Find sharp lower bounds on the probability*

$$\mathbb{P}\{H_*(N^\pm) \cong H_*(Q_M^\pm)\}.$$

While this setting might seem restrictive, it does arise naturally in the context of evolution equations. If we consider linear stochastic partial dif-ferential equations or linear deterministic partial differential equations originating at randomly chosen initial conditions, then the solution at any later point in time will be a random Fourier series, as in (1).

We would also like to point out that our setting is a special case of random fields. The geometry of random fields has been studied by a variety of authors; see, for example, [1, 2] and the references therein. These studies, however, almost exclusively consider the study of the Euler characteristic of nodal domains, which are called *excursion sets* in their setting. As was pointed out in [19], however, in many cases studying the Euler characteristic alone does not provide enough information. We therefore hope that our study of the Betti numbers in this paper will open the door to a more refined study of geometric properties in a stochastic setting.

With these introductory comments in mind, let us now turn to the results of this paper. To focus on obtaining sharp lower bounds asked for in (P), we restrict our attention to one- and two-dimensional domains.

THEOREM 1.4. *Consider a probability space $(\Omega, \mathcal{F}, \mathbb{P})$, a compact interval $G = [a,b]$ in $\mathbb{R}$ and a random field $u: G \times \Omega \to \mathbb{R}$ over $(\Omega, \mathcal{F}, \mathbb{P})$ for which assumptions* (A1) *and* (A2) *hold and such that for $\mathbb{P}$-almost all $\omega \in \Omega$, the function $u(\cdot, \omega): G \to \mathbb{R}$ is twice continuously differentiable. In addition, for $x \in G$ and $\delta > 0$ with $x + \delta \in G$, let*

$$p_\sigma(x, \delta) = \mathbb{P}\left\{\sigma \cdot u(x, \omega) \geq 0, \sigma \cdot u\left(x + \frac{\delta}{2}, \omega\right) \leq 0, \sigma \cdot u(x + \delta, \omega) \geq 0\right\}$$

*and assume that there exists a constant $\mathcal{C}_0 > 0$ such that*

(2) $\quad p_\sigma(x, \delta) \leq \mathcal{C}_0 \cdot \delta^3 \qquad$ *for all $\sigma \in \{\pm 1\}$ and $x \in G$ with $x + \delta \in G$.*

*Then for every discretization size $M$, the probability that the homologies of $N^\pm(\omega)$ and $Q_M^\pm(\omega)$ coincide satisfies*

(3) $$\mathbb{P}\{H_*(N^\pm) \cong H_*(Q_M^\pm)\} \geq 1 - \frac{8\mathcal{C}_0(b-a)^3}{3M^2}.$$

The proof of this theorem—which is fairly straightforward when viewed from the proper perspective—and its implications are presented in Section 2. Observe that the assumption in (2) lies at the heart of Theorem 1.4 and establishing its validity generally requires some work. In particular, determining or estimating the constant $\mathcal{C}_0$ has to be done with significant care since it has direct implications for the tightness of the resulting bounds. We obtain a value for $\mathcal{C}_0$ in the case of *random periodic functions* given by random Fourier series with independent Gaussian coefficients (Theorem 2.7) and indicate the implications for *random trigonometric polynomials*.

The corresponding result for two-dimensional domains, Theorem 3.8, is considerably more subtle to state and left to Section 3. In particular, as is explained in Section 3.1, the obvious generalization of the assumption in (2) leads to a suboptimal choice of exponent for $M$. We remark that the results of Theorem 3.8 depend on two constants $\mathcal{C}_1$ and $\mathcal{C}_2$ associated with the probabilistic behavior of the function at points on the boundary and the interior of $G$, respectively. Values for these constants are computed in the case of *random doubly periodic functions* and applied to the case of random trigonometric polynomials in two variables.

Not surprisingly, our results require several technical probabilistic estimates which are presented in Section 4.

**2. The one-dimensional case.** This section contains the proof of Theorem 1.4, the first step of which is to provide a finite criterion against which one can determine if the homology computation is correct. Using this criterion, a value for $\mathcal{C}_0$ in the case of random periodic functions given by random Fourier series with independent Gaussian coefficients is computed. This result is then applied to the simpler setting of random trigonometric polynomials. The latter theoretical results are then compared against data obtained from numerical simulations.

2.1. *Double crossovers and homology validation.* The first step in addressing the correctness of homology computations for random functions is to formulate a validation criterion which is amenable to a probabilistic study. This will be accomplished in the current subsection, in a deterministic setting.

Let $G = [a,b] \subset \mathbb{R}$ denote a compact interval and let $u : G \to \mathbb{R}$ denote an arbitrary continuous function. We are interested in determining the homology of the nodal domains

$$N^\pm = \{x \in [a,b] : \pm u(x) \geq 0\},$$

that is, we would like to determine the number of components of each nodal domain. As mentioned in the introduction, rather than working with $N^\pm$ directly, we choose a discretization size $M$ and study the associated cubical approximations $Q_M^\pm$ introduced in Definition 1.1. These sets depend only on the function values of $u$ at the $M+1$ discretization points $x_k = a + k(b-a)/M$, where $k = 0, \ldots, M$. For sufficiently large $M$ and under suitable regularity and nondegeneracy conditions on $u$, one can then expect the number of components of $Q_M^\pm$ and $N^\pm$ to agree.

In the one-dimensional situation, of course, it is straightforward to characterize when the homologies of $Q_M^\pm$ coincide with the homologies of the nodal domains $N^\pm$. This is the case if, on each subinterval of the form $[x_k, x_{k+1}]$:

- the function $u$ has no zero at all, or alternatively,

- it has exactly one zero in the interior of the interval and different signs on either side of the zero.

In other words, if on each interval $(x_k, x_{k+1})$ between two consecutive discretization points, the function $u$ has at most one zero and if $u$ is nonzero at all $x_k$, then the homology computation is correct. This shows that in the context of random functions, we have to find sharp upper bounds on the probability of having at least two zeros (counting multiplicities) in a given interval. One would expect that for small intervals, this probability will be small.

Unfortunately, the above characterization in terms of zeros is not directly applicable to random functions since in this context, the condition has to be formulated in terms of the finite-dimensional distributions of the random function. One possible approach is based on the following concept.

DEFINITION 2.1 (*Double crossover*). Let $u:[a,b] \to \mathbb{R}$ be a continuous function and let $[\alpha, \beta] \subset [a,b]$ be arbitrary. Then *u has a double crossover on* $[\alpha, \beta]$ if

$$(4) \qquad \sigma \cdot u(\alpha) \geq 0, \qquad \sigma \cdot u\left(\frac{\alpha+\beta}{2}\right) \leq 0 \quad \text{and} \quad \sigma \cdot u(\beta) \geq 0$$

for one choice of the sign $\sigma \in \{\pm 1\}$.

Since the existence of a double crossover only depends on three function values of the function $u$, the probability for the occurrence of a double crossover can be estimated for random functions. At first glance, however, the obvious connection between double crossovers and zeros is of little use for our problem. Since the existence of a double crossover on $[\alpha, \beta]$ implies the existence of at least two zeros, the double crossover probability can only provide a lower bound on the probability of having at least two zeros in $[\alpha, \beta]$—yet we are interested in an upper bound. In order to accomplish this, we use a device which was introduced in [15].

DEFINITION 2.2 (*Dyadic subintervals, admissible intervals*). Let $u:[a,b] \to \mathbb{R}$ be a continuous function and let $J = [\alpha, \beta] \subset [a,b]$ be arbitrary.

- The *dyadic points* in the interval $J = [\alpha, \beta]$ are defined as

$$d_{n,k} = \alpha + (\beta - \alpha) \cdot \frac{k}{2^n} \qquad \text{for all } k = 0, \ldots, 2^n \text{ and } n \in \mathbb{N}_0.$$

The *dyadic subintervals* of $J$ are the intervals $[d_{n,k}, d_{n,k+1}]$ for all $k = 0, \ldots, 2^n - 1$ and $n \in \mathbb{N}_0$.
- The interval $J = [\alpha, \beta]$ is called *admissible for u* if the function $u$ does not have a double crossover on any of the dyadic subintervals of $J$.

The notion of admissibility depends only on the signs of the function values at the countable set of dyadic points in the interval, which makes it ideal for a probabilistic setting. Moreover, it allows us to determine the locations of zeros of a continuous function. We begin by presenting a criterion which rules out sign changes.

LEMMA 2.3. *Let $u: J \to \mathbb{R}$ be a continuous function on $J = [\alpha, \beta]$ and assume that $J$ is admissible for $u$. Furthermore, suppose that both $u(\alpha) \geq 0$ and $u(\beta) \geq 0$. Then for all $x \in J$, we have $u(x) \geq 0$. Analogously, the inequalities $u(\alpha) \leq 0$ and $u(\beta) \leq 0$ imply that $u(x) \leq 0$ on $J$.*

PROOF. Assume that $u(\alpha) \geq 0$ and $u(\beta) \geq 0$. Since the function $u$ is continuous, it suffices to verify inductively that $u(d_{n,k}) \geq 0$ at all dyadic points $d_{n,k}$. According to our assumption, we have $u(d_{0,k}) \geq 0$ for $k = 0, 1$.

Now, assume that for some $n \in \mathbb{N}_0$, we have shown $u(d_{n,k}) \geq 0$ for all $k = 0, \ldots, 2^n$. Since $J$ is admissible, the function $u$ has no double crossover on the dyadic subinterval $[d_{n,k}, d_{n,k+1}]$ for any $k = 0, \ldots, 2^n - 1$. With $u(d_{n,k}) \geq 0$ and $u(d_{n,k+1}) \geq 0$, this immediately implies that $u((d_{n,k} + d_{n,k+1})/2) \geq 0$ and this in turn implies that $u(d_{n+1,\ell}) \geq 0$ for all $\ell = 0, \ldots, 2^{n+1}$. $\square$

The next lemma addresses the possibility of a sign change on an admissible interval.

LEMMA 2.4. *Let $u: J \to \mathbb{R}$ be a continuous function on $J = [\alpha, \beta]$ and assume that $J$ is admissible for $u$. Furthermore, suppose that $u(\alpha) \geq 0 \geq u(\beta)$. Then there exists a point $x^* \in (\alpha, \beta)$ such that*

$$u(x) \geq 0 \quad \text{for all } \alpha \leq x \leq x^* \quad \text{and} \quad u(x) \leq 0 \quad \text{for all } x^* \leq x \leq \beta.$$

*An analogous result holds if $u(\alpha) \leq 0 \leq u(\beta)$.*

PROOF. We recursively introduce two sequences $(\alpha_k)_{k \in \mathbb{N}_0}$ and $(\beta_k)_{k \in \mathbb{N}_0}$ in $J$ as follows. Let $\alpha_0 = \alpha$, $\beta_0 = \beta$ and for $k \in \mathbb{N}_0$, set

$$\alpha_{k+1} = \alpha_k \quad \text{and} \quad \beta_{k+1} = \frac{\alpha_k + \beta_k}{2} \quad \text{if } u\left(\frac{\alpha_k + \beta_k}{2}\right) \leq 0,$$

$$\alpha_{k+1} = \frac{\alpha_k + \beta_k}{2} \quad \text{and} \quad \beta_{k+1} = \beta_k, \quad \text{if } u\left(\frac{\alpha_k + \beta_k}{2}\right) > 0.$$

Then the sequence $(\alpha_k)$ is increasing and the sequence $(\beta_k)$ is decreasing. In addition, we have both $\beta_k - \alpha_k = (\beta - \alpha)/2^k$ and $u(\alpha_k) \geq 0 \geq u(\beta_k)$ for all $k \in \mathbb{N}_0$. Using Lemma 2.3, it is also straightforward to verify that for all $k \in \mathbb{N}_0$, we have both $u(x) \geq 0$ on $[\alpha, \alpha_k]$ and $u(x) \leq 0$ on $[\beta_k, \beta]$. The result now follows with $x^* = \lim_{k \to \infty} \alpha_k$. $\square$

Note that in the situation of the above lemma, the function $u$ can certainly have more than one zero. In addition to the zero $x^*$, it is possible that $u$ has zeros in both $(\alpha, x^*)$ and $(x^*, \beta)$, provided the function does not change sign at these zeros. In other words, any additional zeros have to be double zeros.

By combining the above two lemmas, we can now formulate our validation criterion for homology computations in one dimension.

PROPOSITION 2.5 (*Validation criterion*). *Let $u: I \to \mathbb{R}$ be a continuous function on the compact interval $I = [a, b]$ and let $M \in \mathbb{N}$ be arbitrary. Let $N^\pm$ denote the nodal domains of $u$ and let $Q_M^\pm$ denote their cubical approximations, as in Definition 1.1. Furthermore, assume that the following hold:*

(a) *the function $u$ is nonzero at all grid points $x_k$ for $k = 0, \ldots, M$;*

(b) *the function $u$ has no double zero in $(a, b)$, that is, if $x \in (a, b)$ is a zero of $u$, then $u$ attains both positive and negative function values in every neighborhood of $x$;*

(c) *for every $k = 0, \ldots, M - 1$, the interval $[x_k, x_{k+1}]$ between consecutive discretization points is admissible for $u$ in the sense of Definition 2.2.*

*We then have*

$$H_*(N^\pm) \cong H_*(Q_M^\pm),$$

*that is, the homologies of the nodal domains and of their cubical approximations coincide.*

PROOF. Consider the interval $[x_k, x_{k+1}]$. If $u(x_k)$ and $u(x_{k+1})$ are both positive, then by Lemma 2.3 and (b), the function $u$ must be strictly positive on $[x_k, x_{k+1}]$. Similarly, negative function values force $u$ to be negative on $[x_k, x_{k+1}]$. Finally, if the function values at the endpoints have opposite signs, then Lemma 2.4 and (b) imply the existence of a unique zero $x^*$ in $(x_k, x_{k+1})$. From this, the result follows easily. □

Proposition 2.5 leads to a straightforward proof of Theorem 1.4.

PROOF OF THEOREM 1.4. Fix a discretization size $M \in \mathbb{N}$ and for $\ell \in \{0, \ldots, M-1\}$, consider the interval $I_\ell = [x_\ell, x_{\ell+1}]$ between two consecutive discretization points, as defined in Definition 1.1. The interval $I_\ell$ is not admissible if the function $u$ has a double crossover on one of its dyadic subintervals. If we now denote the dyadic points in $I_\ell$ by $d_{n,k}$, as in Definition 2.2, then together with (2) and $\delta = (b-a)/M$, we obtain the estimate

$$\mathbb{P}\{I_\ell \text{ is not admissible}\} \leq \sum_{n=0}^{\infty} \sum_{k=0}^{2^n - 1} \mathbb{P}\{u \text{ has a double crossover in } [d_{n,k}, d_{n,k+1}]\}$$

$$\leq \sum_{n=0}^{\infty} \sum_{k=0}^{2^n-1} (p_{+1}(d_{n,k}, \delta/2^n) + p_{-1}(d_{n,k}, \delta/2^n))$$

$$\leq \sum_{n=0}^{\infty} \sum_{k=0}^{2^n-1} 2\mathcal{C}_0 \cdot \left(\frac{\delta}{2^n}\right)^3 = \frac{8\mathcal{C}_0(b-a)^3}{3M^3}.$$

According to Proposition 2.5 and Assumptions (A1) and (A2), the probability that the homologies of $N^{\pm}(\omega)$ and $Q_M^{\pm}(\omega)$ differ can then be estimated as

$$1 - \mathbb{P}\{H_*(N^{\pm}) \cong H_*(Q_M^{\pm})\} \leq \mathbb{P}\{0 \text{ is a critical value}\} + \sum_{\ell=0}^{M} \mathbb{P}\{u(x_\ell) = 0\}$$

$$+ \sum_{\ell=0}^{M-1} \mathbb{P}\{I_\ell \text{ is not admissible}\}$$

$$\leq \frac{8\mathcal{C}_0(b-a)^3}{3M^2},$$

which completes the proof of the theorem. □

Note that the above proof could easily be adapted to also cover nonuniform discretizations of the domain $G$.

2.2. *Random periodic functions.* In order to apply Theorem 1.4 to specific random sums, we must verify the assumption in (2). For the case of periodic random functions, this is demonstrated in the following.

ASSUMPTION 2.6. Let $(\Omega, \mathcal{F}, \mathbb{P})$ denote a probability space, define $G = [0, L]$ for some fixed $L > 0$ and consider the random Fourier series $u : G \times \Omega \to \mathbb{R}$ defined by

$$(5) \qquad u(x,\omega) = \sum_{k=0}^{\infty} a_k \cdot \left(g_{2k}(\omega) \cdot \cos\frac{2\pi k x}{L} + g_{2k-1}(\omega) \cdot \sin\frac{2\pi k x}{L}\right).$$

We assume that the following hold:

(a) the Gaussian random variables $g_\ell : \Omega \to \mathbb{R}$, $\ell \in \mathbb{N}_0$, in (5) are defined over the common probability space $(\Omega, \mathcal{F}, \mathbb{P})$ and are independent and normally distributed with mean 0 and variance 1;

(b) the constants $a_k$, $k \in \mathbb{N}_0$, in (5) are arbitrary real numbers such that at least two of them do not vanish and such that

$$\sum_{k=0}^{\infty} k^6 a_k^2 < \infty.$$

We would like to point out that in the above assumption, the constants $a_k$ are directly related to smoothness properties of the random function $u$. More precisely, one can show that if

$$\sum_{k=0}^{\infty} k^{2p} a_k^2 < \infty \quad \text{for some } p > 0,$$

then $\mathbb{P}$-almost all realizations $u(\cdot, \omega)$ are contained in the Hölder space $C^q[0, L]$ for any real $0 < q < p$; see, for example, [22], Section 7.4. Moreover, if we define

(6) $$A_\ell = \sum_{k=0}^{\infty} k^{2\ell} a_k^2,$$

then it is straightforward to verify that, in fact,

$$\mathbb{E} \|D_x^\ell u\|_{L^2(0,L)}^2 = (2\pi)^{2\ell} \cdot L^{1-2\ell} \cdot A_\ell,$$

where $\mathbb{E}$ denotes the expected value of a random variable over $(\Omega, \mathcal{F}, \mathbb{P})$. In other words, the constant $A_\ell$ contains averaged information on the $L^2(0, L)$-norm of the $\ell$th derivative of the random function $u$.

Random periodic functions as in (5) are particularly convenient to study due to special properties of their spatial covariance functions. To see this, consider the spatial covariance function of (5), which is given by

(7) $$R(x, y) = \mathbb{E} u(x) u(y) = \sum_{k=0}^{\infty} a_k^2 \cdot \cos \frac{2\pi k(x - y)}{L}.$$

In other words, under the above assumptions, the random function defined in (5) is a homogeneous random field in the sense of [1]—this simplifies some of the considerations significantly.

Now consider the problem of determining the homology of the nodal domains of $u$ using a spatial discretization of size $M$.

THEOREM 2.7 (Random periodic functions in 1D). *We consider the random $L$-periodic function $u : G \times \Omega \to \mathbb{R}$ defined in (5), where $G = [0, L]$ with $L > 0$, and suppose that Assumption 2.6 holds.*

*Let $M$ denote an arbitrary positive integer and, as in Definition 1.1, let $Q_M^\pm(\omega)$ denote the cubical approximations of the random nodal domains $N^\pm(\omega)$ of $u(\cdot, \omega)$ which are used for the computation of the homologies of $N^\pm(\omega)$. Then the probability that the homology of the random nodal domains $N^\pm(\omega)$ is computed correctly with the discretization of size $M$ satisfies*

(8) $$\mathbb{P}\{H_*(N^\pm) \cong H_*(Q_M^\pm)\} \geq 1 - \frac{\pi^2}{6M^2} \cdot \frac{A_0 A_2 - A_1^2}{A_0^{3/2} A_1^{1/2}} + O\left(\frac{1}{M^3}\right),$$

*where the constants $A_\ell$ are defined as in (6).*

PROOF. According to our assumptions, the random variable $u(x, \cdot) : \Omega \to \mathbb{R}$ is normally distributed with mean 0 and variance $\sum_{k=0}^{\infty} a_k^2 \neq 0$ for each $x \in G$, which immediately implies (A1). Furthermore, (A2) follows readily from [1], Theorem 3.2.1. Thus in order to apply Theorem 1.4, we need only verify (2).

We employ Proposition 4.1 with $n = 3$ and sign vector $(s_1, s_2, s_3) = (1, -1, 1)$. Let $x \in G$ be arbitrary, but fixed, and consider the $\delta$-dependent three-dimensional random vector $T(\delta) = (T_1(\delta), T_2(\delta), T_3(\delta)) : \Omega \to \mathbb{R}^3$ defined by

$$T_1(\delta) = u(x), \qquad T_2(\delta) = u\left(x + \frac{\delta}{2}\right) \quad \text{and} \quad T_3(\delta) = u(x + \delta).$$

Then $T(\delta)$ is a Gaussian random variable with mean $0 \in \mathbb{R}^3$ and covariance matrix $C(\delta)$ given by

$$C(\delta) = \begin{bmatrix} r(0) & r(\delta/2) & r(\delta) \\ r(\delta/2) & r(0) & r(\delta/2) \\ r(\delta) & r(\delta/2) & r(0) \end{bmatrix},$$

where, in view of (7), we use the abbreviation

$$r(\delta) = R(y + \delta, y) = \sum_{k=0}^{\infty} a_k^2 \cdot \cos \frac{2\pi k \delta}{L}.$$

For $\delta \to 0$, the even function $r$ can be expanded as

$$r(\delta) = R_0 - \frac{R_1}{2} \cdot \delta^2 + \frac{R_2}{24} \cdot \delta^4 - \frac{R_3}{720} \cdot \delta^6 + O(\delta^8) \qquad \text{with } R_\ell = \left(\frac{2\pi}{L}\right)^{2\ell} \cdot A_\ell$$

and Assumption 2.6(b) implies that

(9) $$R_0 \neq 0, \qquad R_1 \neq 0 \quad \text{and} \quad R_0 R_2 - R_1^2 > 0.$$

For the last inequality, one need only observe that the Cauchy–Schwarz inequality implies $R_1^2 \leq R_0 R_2$, and that equality can easily be excluded since two of the $a_k$ are nonzero. Using the above expansion, the determinant of the covariance matrix can be written as

(10) $$\det C(\delta) = \frac{R_1}{64} \cdot (R_0 R_2 - R_1^2) \cdot \delta^6 + O(\delta^8),$$

that is, the covariance matrix is positive definite for sufficiently small $\delta > 0$ and Proposition 4.1(i) is satisfied. Furthermore, by applying the Newton polygon method [27, 29] to the characteristic polynomial $\det(C(\delta) - \lambda I)$, it can be shown that in the limit $\delta \to 0$, the three eigenvalues $\lambda_k(\delta)$, for $k = 1, 2, 3$, of $C(\delta)$ satisfy

$$\lambda_1(\delta) = \frac{R_0 R_2 - R_1^2}{96 R_0} \cdot \delta^4 + O(\delta^5),$$

$$\text{(11)} \quad \lambda_2(\delta) = \frac{R_1}{2} \cdot \delta^2 + O(\delta^3),$$
$$\lambda_3(\delta) = 3R_0 + O(\delta).$$

This establishes the validity of Proposition 4.1(iii). In order to determine the asymptotic behavior of the eigenvector corresponding to $\lambda_1$, we consider the adjoint of the covariance matrix, whose expansion is given by

$$\operatorname{adj} C(\delta) = \frac{R_0 R_1}{4} \cdot \begin{bmatrix} 1 & -2 & 1 \\ -2 & 4 & -2 \\ 1 & -2 & 1 \end{bmatrix} \cdot \delta^2 + O(\delta^3).$$

The constant coefficient matrix has the double eigenvalue 0, as well as the positive eigenvalue 6 with associated unnormalized eigenvector $(1, -2, 1)^t$. Since the eigenspace for the largest eigenvalue of the adjoint matrix coincides with the eigenspace for the eigenvalue $\lambda_1(\delta)$ of $C(\delta)$, the simplicity of these eigenvalues shows that we can choose an eigenvector $v_1(\delta)$ for $\lambda_1(\delta)$ which satisfies $v_1(\delta) \to (1, -2, 1)^t / 6^{1/2}$ for $\delta \to 0$ and therefore Proposition 4.1(ii) is satisfied. Applying Proposition 4.1, we finally obtain

$$\lim_{\delta \to 0} p_{+1}(x, \delta) \cdot \sqrt{\frac{\det C(\delta)}{\lambda_1(\delta)^3}} = \frac{3\sqrt{6}}{8\pi}.$$

In combination with (9), (10) and (11), this limit implies that

$$p_\sigma(x, \delta) = \frac{1}{128\pi} \cdot \frac{R_0 R_2 - R_1^2}{R_0^{3/2} R_1^{1/2}} \cdot \delta^3 + O(\delta^4) = \frac{\pi^2}{16L^3} \cdot \frac{A_0 A_2 - A_1^2}{A_0^{3/2} A_1^{1/2}} \cdot \delta^3 + O(\delta^4)$$

for $\sigma = 1$, the analogous statement for $\sigma = -1$ following from a symmetry argument. Thus (2) is also satisfied and Theorem 2.7 follows from our abstract result, Theorem 1.4. □

Very much in the spirit of [25], our result furnishes an explicit bound on the likelihood of computing the correct homology. The bound depends on the discretization size $M$ and global properties of the underlying function $u$. Unlike [25], however, the necessary information on $u$ can be easily computed. In fact, it is given explicitly in terms of the coefficient sequence $(a_k)$ which, in turn, is related to smoothness properties of $u$. Note also that due to the asymptotic nature of our results, the estimate in (8) contains an additional $O(1/M^3)$ term. While in principle it is also possible to derive explicit bounds for this term, the formulas involved are rather complicated. In practice, however, the term can already be neglected for relatively small discretization sizes. For more details, we refer the reader to [11, 12].

2.3. *Random trigonometric polynomials.* We now demonstrate the implications of our main results for the case of random trigonometric polynomials. Let $N \geq 3$ denote an arbitrary integer and let $L = 2\pi$. If, in the situation of Theorem 2.7, we choose the coefficients specifically as $a_k = 1$ for $1 \leq k \leq N$, $a_0 = 0$ and $a_k = 0$ for $k > N$, then (5) represents a random trigonometric polynomial,

$$(12) \qquad u(x,\omega) = \sum_{k=1}^{N} (g_{2k}(\omega) \cdot \cos(kx) + g_{2k-1}(\omega) \cdot \sin(kx)).$$

The asymptotic behavior of the number of zeros of such functions has been studied extensively and therefore provides us with a means for testing Theorem 2.7. Specifically, one can show that as $N \to \infty$, most trigonometric polynomials of the form (12) have on the order of $2N/3^{1/2}$ real zeros in the interval $[0, 2\pi]$; see, for example, [4, 16, 17]. Thus in order for our homology computation to be correct, we must at least ensure that the discretization size $M$ satisfies $1/M = O(1/N)$. Due to the almost certain occurrence of zeros which are more closely spaced, it seems implausible that this bound would be optimal. Now, let $A_\ell$ be defined as in (6). Then one can easily show that

$$\frac{A_0 A_2 - A_1^2}{A_0^{3/2} A_1^{1/2}} = \frac{\sqrt{6}}{180}(N-1)(8N+11)\sqrt{(N+1)(2N+1)} \sim \frac{4\sqrt{3}}{45} \cdot N^3.$$

This implies that in order for the homology computation to be accurate with high confidence, we must choose the discretization size $M$ in such a way that

$$M \sim N^{3/2} \qquad \text{for } N \to \infty.$$

We conjecture that this proportionality is the correct asymptotic result, as is supported by the numerical results shown in Figure 4. For various values of $N$ between 5 and 1000, we computed three different quantities: the expected number of zeros (lowermost curve), the expected minimal distance between two consecutive zeros of the random polynomial (second curve from bottom) and the integer $M$ for which 95% of the generated random polynomials have a minimal distance between consecutive zeros which is larger that $2\pi/M$ (uppermost curve). For comparison, the prediction of Theorem 2.7 which would ensure a correctness probability of at least 95% is given by the curve, which lies exactly below the uppermost curve and shows the same asymptotic growth behavior.

Numerically validated computational results which also treat the case of linear and nonlinear evolution equations (such as the Cahn–Hilliard models described in the introduction) will be presented in [12].

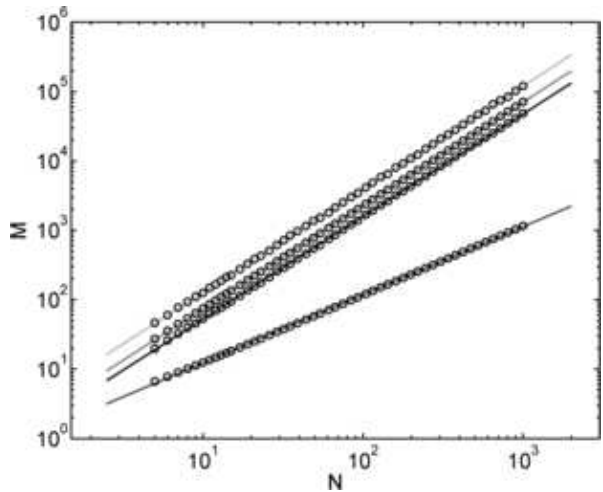

Fig. 4. *Numerical results for random trigonometric polynomials.*

**3. The two-dimensional case.** In this section, we address the case of two-dimensional base domains $G \subset \mathbb{R}^2$ and are interested in the topology of the nodal domains of a function $u : G \to \mathbb{R}$. For simplicity, it is assumed that $G$ is a square domain, even though more general domains could also be treated. As in the one-dimensional situation, we present both general results and specific applications, in particular, the case of doubly periodic random functions.

3.1. *Double crossovers and $B$-admissibility.* In contrast to the one-dimensional situation where a nodal domain of a function consists of a disjoint union of intervals, the geometry of nodal domains can be rather complicated in two space dimensions. It is therefore not surprising that the method of Section 2 cannot be generalized directly to the higher-dimensional case. In this first subsection, we show that a fairly straightforward adaptation of the one-dimensional ideas can be used to discuss the accuracy of homology computations. Unfortunately, however, this approach alone will provide suboptimal results and must therefore be extended significantly in the following subsection.

In the one-dimensional situation, the basic ingredients of our approach can be summarized as follows. After recognizing the fundamental obstruction to a valid homology computation as the occurrence of two zeros in an interval between gridpoints, this obstruction was reformulated in terms of the concept of a double crossover. The latter is a specific function value sign pattern that can be observed at the endpoints and midpoint of a given interval $J$. Moreover, by excluding this sign pattern from all dyadic subintervals of $J$, we were able to introduce the notion of admissibility of $J$. Admissibility guarantees that the topological information for the nodal domains in $J$ can be determined from the function values of the given function at the endpoints of $J$. If the function values at the endpoints have the same sign, then the function has no zero in $J$; if the signs differ, then the function has exactly one zero in $J$.

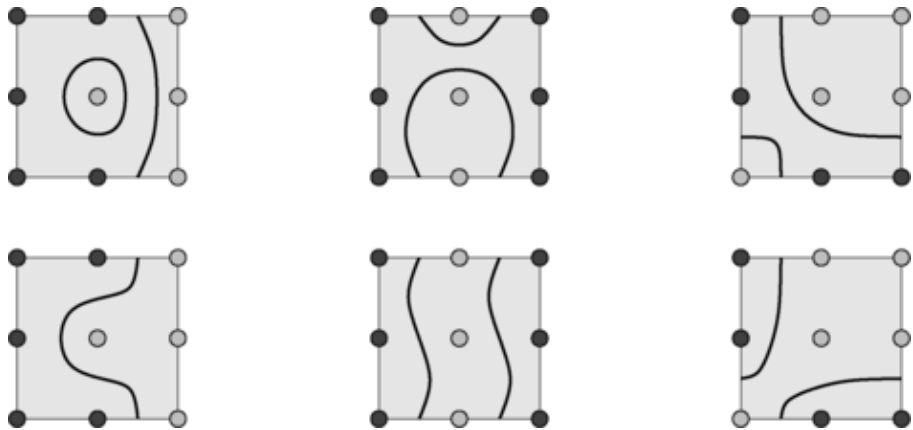

FIG. 5. *Examples of problem situations in two dimensions. In each of the three columns, the function values at the grid points shown have identical sign configurations, but the underlying topology of the nodal domains differs. The nodal lines are shown as solid black lines.*

Based on the above discussion, our goal is the identification of certain sign patterns of function values that could indicate mistakes in the homology computation. As the fundamental building block, we consider a square with side length $\delta$ and, analogously to the one-dimensional case, we consider the function values of the underlying function $u$ at the corners as well as at the midpoints of the sides and at the center of the square. This will easily enable us to consider dyadic subsquares as before.

Some representative examples of problematic sign patterns are shown in Figure 5. Within each column, the two examples exhibit the same sign pattern, but the underlying topology of the nodal domains differs. The examples in the first two columns are, in some sense, reminiscent of the one-dimensional situation since in each of these cases, we can find a double crossover on at least one vertical or horizontal line of grid points. This is not true for the examples in the last column of Figure 5. Here, there are no double crossovers on vertical or horizontal lines, yet the underlying topology still differs. On the level of the function values of $u$, this is due to the fact that, at the corners, the signs alternate. In other words, there are two double crossovers around the perimeter of the square. Altogether, we see that each of the sign configurations shown in Figure 5 contains one of the seven sign patterns depicted in Figure 6.

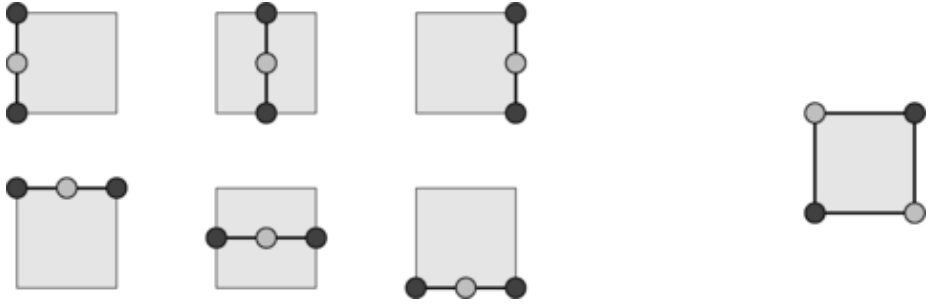

FIG. 6. *For the notion of B-admissibility, the seven sign configurations shown above are forbidden. Each image corresponds to two sign patterns: one in which the dark points correspond to positive function values and the light points correspond to negative ones, and one for the opposite situation.*

It turns out that excluding the patterns in Figure 6 is sufficient to extend the general procedure of Section 2. For this, we need the following definitions.

DEFINITION 3.1 (*Dyadic subsquares*). Let $J := [\alpha, \alpha+\delta] \times [\beta, \beta+\delta] \subset \mathbb{R}^2$ and denote the dyadic points in the interval $[0,\delta]$ in the sense of Definition 2.2, by $d_{n,k}$. Then the *dyadic points* in the square $J$ are the points

$$d_{n,k,\ell} = (\alpha + d_{n,k}, \beta + d_{n,\ell}) \in \mathbb{R}^2 \qquad \text{for all } k, \ell = 0, \ldots, 2^n \text{ and } n \in \mathbb{N}_0.$$

The *dyadic subsquares of* $J$ are the squares $d_{n,k,\ell} + [0, \delta/2^n]^2$ for all $k, \ell = 0, \ldots, 2^n - 1$ and $n \in \mathbb{N}_0$.

DEFINITION 3.2 (*B-admissibility*). Let $u: G \to \mathbb{R}$ denote an arbitrary function and let $J \subset G$ denote a square, as in Definition 3.1. Then $J$ is called *B-admissible for* $u$, if none of its dyadic subsquares contain any of the function value sign configurations shown in Figure 6.

Using the concept of $B$-admissibility, one can validate homology more or less analogously to the one-dimensional case. To demonstrate this, consider Figure 7. If we assume that a given square $J$ of side length $\delta > 0$ (shown in pink) is $B$-admissible for $u$, then the figure shows all possible function value sign configurations which can be observed at the dyadic points $d_{1,k,\ell}$, where $k, \ell = 0, 1, 2$, up to symmetry actions of the square. These configurations are grouped according to the sign configurations at the corners of the square. We now make two important observations:

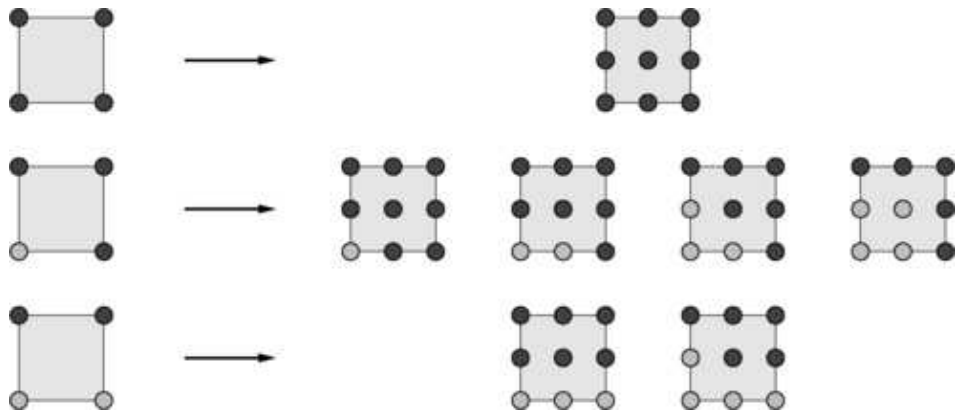

FIG. 7. *Possible sign patterns for the function values of $u$ at the dyadic points $d_{1,k,\ell}$, where $k,\ell = 0,1,2$, if the underlying square is $B$-admissible. In addition to the patterns shown, all sign patterns which can be obtained by applying an element of the symmetry group of the square can also occur, resulting in a total of 66 possible patterns.*

- if all function values of $u$ at the four corners of $J$ have the same sign, then all function values at the nine depicted dyadic points must necessarily have the same sign;
- if both signs can be observed at the corners of $J$, then there are exactly two sides of the square with contain both positive and negative function values.

The first observation is reminiscent of Lemma 2.3. It shows that if $J$ is $B$-admissible for $u$ and $u$ is positive at the corners of $J$, then $u$ cannot take negative function values on $J$. Thus the first observation can be used to identify subsets of one or the other nodal domain. In combination with the second observation, it therefore seems plausible to expect that the nodal line in a $B$-admissible square can be pinned down. In fact, under suitable nondegeneracy assumptions on $u$, it is possible to prove the following result.

PROPOSITION 3.3. *Let $G \subset \mathbb{R}^2$ denote a square, let $u: G \to \mathbb{R}$ be $C^2$ such that $0$ is not a critical value of $u$ and let $J \subset G$ be a $B$-admissible square for $u$. Then the following hold.*

(a) *If $u$ is strictly positive at the corners of $J$, then $u$ is strictly positive on $J$. Similarly, if $u$ is strictly negative at the corners of $J$, then $u$ is strictly negative on $J$.*

(b) *If $u$ takes both positive and negative function values at the corners of $J$, then the nodal line of $u$ in $J$ is a simple smooth curve which connects one side of $J$ with another side of $J$.*

*Now, assume the situation of Definition 1.1 with a positive integer $M$. If all $M^2$ subsquares which are generated by adjacent gridpoints of the equidistant $M$-discretization of $G$ are $B$-admissible for $u$, then we have*

$$H_*(N^\pm) \cong H_*(Q_M^\pm),$$

*that is, the homologies of the nodal domains $N^\pm$ and of their cubical approximations $Q_M^\pm$ coincide.*

The proof of the above result can easily be formulated using ideas similar to those for the one-dimensional case and is therefore omitted.

At first glance, Proposition 3.3 is exactly the validation criterion which allows us to consider the two-dimensional random field case. It is, in fact, possible to proceed as in Section 2 and formulate a probabilistic estimate based on the above proposition. Unfortunately, the resulting theorem provides highly suboptimal probability estimates. This suboptimality stems from the fact that the three-point sign configurations in Figure 6 are not sufficiently unlikely. It turns out that in many applications, each of these configurations has probability of order $O(\delta^3)$, where $\delta$ denotes the side length of the square $J$. On the other hand, the four-point configuration in Figure 6 has probability of order $O(\delta^4)$. Together, the probability for $B$-admissibility is of order $O(\delta^3)$. For the homology validation, we must ensure that each of the $M^2$ subsquares of $G$ are $B$-admissible and that $\delta \sim 1/M$. This finally gives a probability estimate of order $O(1/M)$, in contrast to the $O(1/M^2)$ estimate in one dimension.

3.2. *Excluded sign patterns and $I$-admissibility.* The discussion of the last paragraph indicates two things. On one hand, one has to avoid the one-dimensional double crossover concept in higher dimensions since its occurrence is too likely. On the other hand, the four-point configuration of Figure 6 does exhibit the correct probability estimate. If we were able to define a notion of admissibility which has probability of order $O(\delta^4)$, then the resulting homology validation result would exhibit the desired quadratic order in $M$.

As a first step toward a suitable admissibility condition, we return to the four-point configuration of Figure 6. It was mentioned above that in many applications, the probability of this sign pattern turns out to be of order $O(\delta^4)$—one would expect that nondegenerate affine transformations of this pattern are similarly rare. This leads naturally to the following notion.

DEFINITION 3.4 ($I_4$-*admissibility*).  Let $u : G \to \mathbb{R}$ and $J = [\alpha, \alpha+\delta] \times [\beta, \beta+\delta] \subset G$. Then $J$ is called $I_4$-*admissible for $u$* if it does not contain any of the function value sign configurations shown in Figure 8.

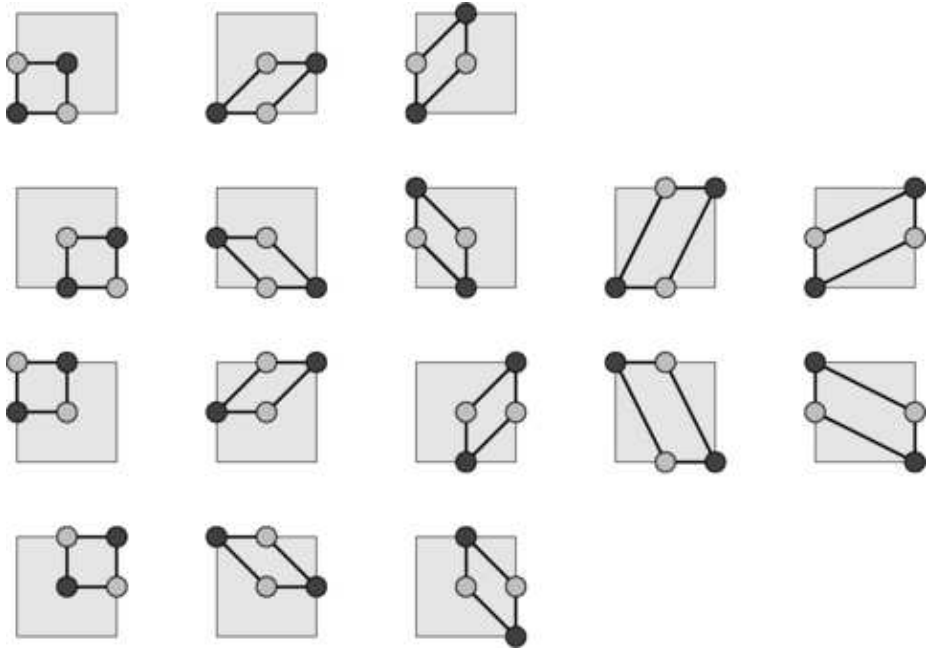

Fig. 8. *For the notion of $I_4$-admissibility, the sixteen sign configurations shown above are forbidden. As in Figure 6, each image corresponds to two sign patterns.*

We are still interested in a recursive approach to homology validation, as described in Sections 2 and 3.1. It is therefore natural to determine which of the 512 sign patterns on the dyadic points $d_{1,k,\ell}$, where $k,\ell = 0,1,2$, can occur if the underlying square is $I_4$-admissible. An exhaustive search shows that only 92 patterns remain possible. Up to symmetry actions of the square, these configurations are shown in Figure 9.

In contrast to the situation described in Section 3.1, not all patterns that are possible under $I_4$-admissibility permit homology validation. The prime example can be seen in the first row of Figure 9. Even though the signs at the corners of the square are all the same, after one level of refinement, the opposite sign appears at the center of the square. This example shows that the concept of $I_4$-admissibility alone cannot be used for homology validation, that is, we must exclude more subconfigurations. At the same time, we must ensure that the probability estimates for the additionally excluded configurations are sufficiently small. As we will see later, it can be shown that in the context of random periodic functions, the last sign configuration in the first row of Figure 9 exhibits an $O(\delta^4)$-bound. In fact, we need not consider the full nine-point pattern, but can rather consider the five-point subpattern shown in Figure 10. This leads to the following definition.

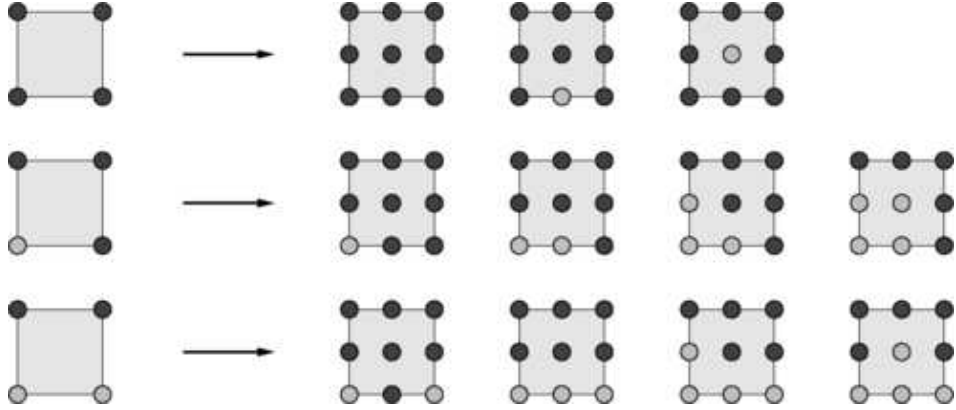

FIG. 9. *Possible sign patterns for the function values of u at the dyadic points $d_{1,k,\ell}$, where $k, \ell = 0, 1, 2$, if the underlying square is $I_4$-admissible. In addition to the patterns shown, all sign patterns which can be obtained by applying an element of the symmetry group of the square can also occur, resulting in a total of 92 possible patterns.*

DEFINITION 3.5 ($I_5$-*admissibility*). Let $u : G \to \mathbb{R}$ and $J = [\alpha, \alpha + \delta] \times [\beta, \beta + \delta] \subset G$. Then $J$ is called $I_5$-*admissible for* $u$ if it does not contain the function value sign configuration shown in Figure 10.

While combining the notions of $I_4$- and $I_5$-admissibility does bring us closer to homology validation, problematic configurations can still be identified. This is demonstrated in Figure 11 by considering groups of two adjacent $I_4$- and $I_5$-admissible squares. In each case, the sign patterns at the corners of the original squares do not suffice to determine the topology after one refinement level.

The topological changes shown in Figure 11 are introduced by specific patterns in Figure 9, namely the second nine-point pattern in the first row and the first nine-point pattern in the last row. Again in the context of random periodic functions, one can determine their likelihood, which provides only the insufficient $O(\delta^3)$-bound. It is therefore impossible to exclude these two patterns directly, or any of their subpatterns. Instead, we follow a different strategy. Note that the two rightmost patterns in Figure 11 can be avoided if the square with center at the midpoint of the common edge between the two original squares is both $I_4$- and $I_5$-admissible. Combining this with the concept of dyadic subsquares (as in Section 3.1) leads to the following, final, definition.

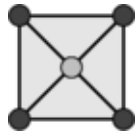

FIG. 10. *For the notion of $I_5$-admissibility, the five-point sign configuration shown above is forbidden. As in Figure 6, this image corresponds to two sign patterns.*

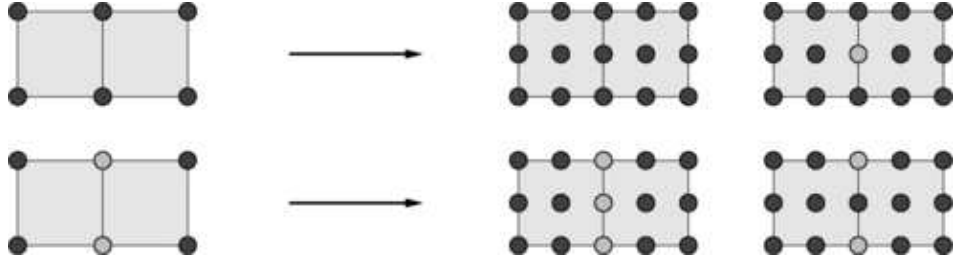

FIG. 11. *Sample configurations which cannot be detected by using both $I_4$- and $I_5$-admissibility, but which exhibit different topological properties.*

DEFINITION 3.6 (*I-admissibility*). Consider $u : G \to \mathbb{R}$ and $J = [\alpha, \alpha + \delta] \times [\beta, \beta + \delta]$ such that the closed $\delta/2$-neighborhood of $J$ is contained in $G$. Then the square $J$ is called *I-admissible for $u$* if for every dyadic subsquare $J^*$ of $J$ with side length $\delta^* = \delta/2^n$, $n \in \mathbb{N}_0$, the following five squares are both $I_4$-admissible and $I_5$-admissible for $u$:

- the dyadic subsquare $J^*$ itself;
- the four shifted squares which are obtained by translating $J^*$ horizontally or vertically by $\delta^*/2$ in either direction, as shown in Figure 12.

As mentioned before, it will be shown later that in many situations, it is possible to derive an $O(\delta^4)$-bound for $I$-admissibility. Note, however, that this improved bound comes at an additional cost. The notion of $I$-admissibility is not determined solely from the function values of $u$ at the dyadic points in $J$. Rather, one must also include knowledge from certain dyadic points in neighboring squares and this renders $I$-admissibility impractical at the boundary of $G$. In order to resolve this final issue, we must combine both $B$- and $I$-admissibility in a suitable way. This is presented in the next subsection.

3.3. *A criterion for homology validation.* In the current subsection, the notions of $B$- and $I$-admissibility are combined to formulate and prove a homology validation criterion in two space dimensions which is amenable to a probabilistic treatment.

While Theorem 1.2 provides a criterion under which it can be guaranteed that for sufficiently large $M$, the homology of the nodal domains agrees with that of the cubical approximation, it does not provide insight into how we would determine a minimum value of $M$. For this, we must consider the notions of $B$- and $I$-admissibility.

PROPOSITION 3.7 (Validation criterion). *Let $G \subset \mathbb{R}^2$ denote a square domain with sides parallel to the coordinate axes and let $u: G \to \mathbb{R}$ be a twice continuously differentiable function for which 0 is a regular value. Let $N^{\pm}$ denote the nodal domains of $u$ and let $Q_M^{\pm}$ denote their cubical approximations (as in Definition 1.1) for some fixed $M \in \mathbb{N}$. In addition, suppose that for each of the $M^2$ closed subsquares $J$ formed by adjacent gridpoints in $G$, $u$ is nonzero at all dyadic subpoints, and that one of the following holds:*

(a) *if the square $J$ lies at the boundary of $G$, then $J$ is B-admissible for $u$;*
(b) *if the square $J$ lies in the interior of $G$, then $J$ is I-admissible for $u$.*

*We then have*

$$H_*(N^{\pm}) \cong H_*(Q_M^{\pm}),$$

*that is, the homologies of the nodal domains and of their cubical approximations coincide.*

PROOF. Let $M$ be fixed as in the formulation of the proposition. In view of Theorem 1.2, it suffices to show that for all values of $n \in \mathbb{N}$, the homologies of $Q_{M \cdot 2^n}^{\pm}$ and $Q_M^{\pm}$ coincide since then we need only choose $n$ sufficiently large for the lemma to hold. In order to establish this invariance of the homologies inductively, we only consider the homologies of $Q_{2M}^{\pm}$ and $Q_M^{\pm}$. Subsequent inductive steps can be carried out completely analogously.

Now consider, without loss of generality, a connected component $K$ of the set $Q_M^+$. We can associate with $K$ a collection of positive grid points $P$ in the underlying discretization of $G$, that is, grid points at which the function $u$ takes positive function values. We must show that the insertion of intermediate grid points in the underlying grid does not change the topology of the component $K$ in the following sense: The corresponding component in $Q_{2M}^+$ has not merged with another component and the number of holes in $K$ has remained constant.

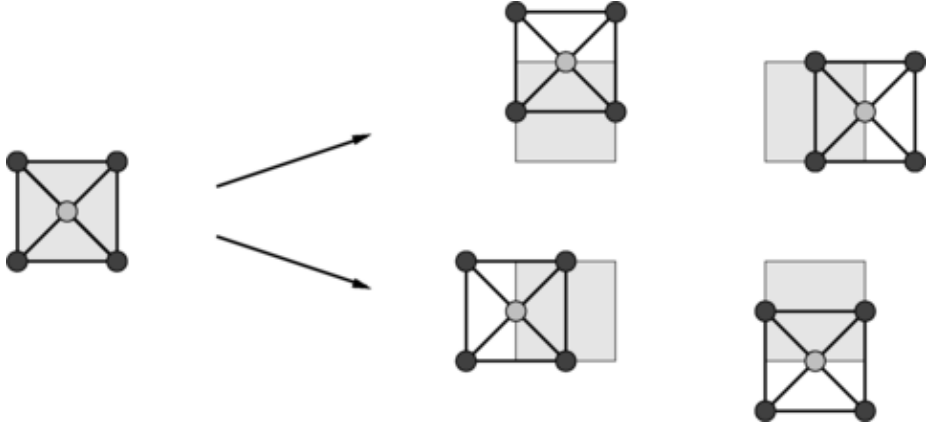

FIG. 12. *For the notion of I-admissibility of a square $J$, we require that every dyadic subsquare $J^*$ (shown with shaded background) of $J$, as well as the four shifted versions of $J^*$ shown on the right, are both $I_4$- and $I_5$-admissible. In the figure, the location of the shifted squares is indicated by the five-point configuration of Figure 10. Nevertheless, all configurations in Figure 8 must be excluded from the five indicated squares.*

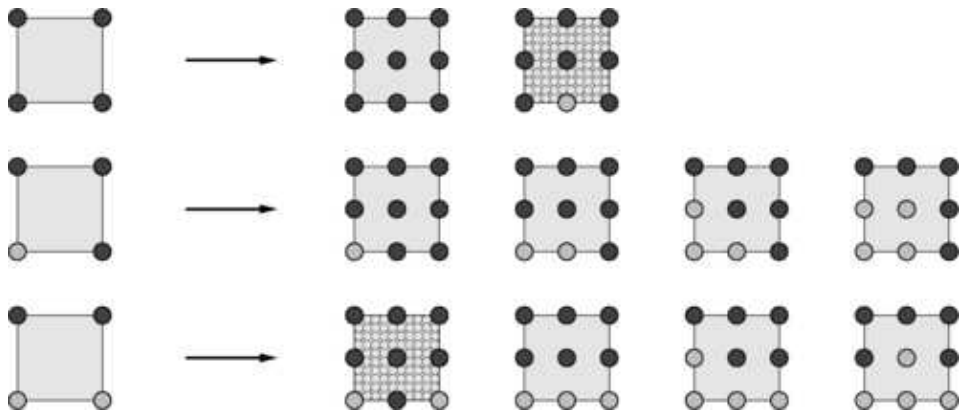

FIG. 13. *Possible sign patterns for the function values of $u$ at the dyadic points $d_{1,k,\ell}$, where $k, \ell = 0, 1, 2$, if the underlying square is $I$-admissible. In addition to the patterns shown, all sign patterns which can be obtained by applying an element of the symmetry group of the square can also occur, resulting in a total of 90 possible patterns. However, the two configurations with patterned backgrounds can only occur in certain situations which are described in more detail in Figure 14.*

In order to perform the step from the given discretization of size $M$ to the one of size $2M$, we must establish which nine-point sign configurations are possible for every given four-point sign configuration at the corners of a discretization square. For the concept of $B$-admissibility, this has already been accomplished in Figure 7. The resulting patterns immediately show that for a $B$-admissible square, passing from $M$ to $2M$ cannot change the topological properties since no new connections are formed and positive/negative grid points are only introduced adjacent to positive/negative components.

As for the concept of $I$-admissibility, we must be more careful. In Figure 13, we collect all possible nine-point sign configurations which could be observed at the next refinement level. With the exception of the two configurations with patterned backgrounds, the same comments as in the last paragraph apply to these patterns and they therefore cannot introduce a change in topology.

Now, consider the two remaining nine-point configurations in Figure 13, that is, the configurations with patterned backgrounds. Recall that our definition of $I$-admissibility includes information about neighboring squares and this implies that these two configurations can only occur if at the edge of the square with the newly generated double crossover, the adjacent square exhibits the specific sign pattern shown in Figure 14. Thus the newly created double crossovers can only occur at the boundary of the component $K$ and if they appear, they do not create any new connections or introduce any change in topology. This fact also implies that in the interior of $K$, no negative gridpoints can be introduced. In addition, note that if the adjacent square is, in fact, a $B$-admissible square, then the two patterned configurations cannot occur at all due to the generation of a double crossover. These facts imply that the topology of the component $K$ will not be changed by passing from $Q_M^+$ to $Q_{2M}^+$ and the proof of the proposition is complete. $\square$

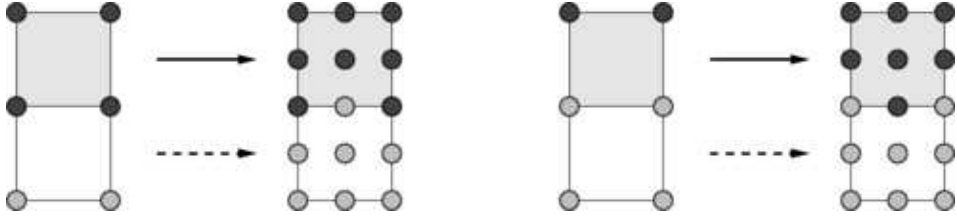

FIG. 14. *The two configurations in Figure 13 with patterned backgrounds can only occur if the neighboring square exhibits a specific sign configuration. For the pattern in the first row of Figure 13, this is shown in the left-hand diagram; for the configuration in the third row see the right-hand diagram.*

We are finally in a position to present an abstract version of a probabilistic homology validation result for two-dimensional random fields $u$ on a square domain $G \subset \mathbb{R}^2$.

THEOREM 3.8 (Abstract 2D estimate). *Consider a probability space $(\Omega, \mathcal{F}, \mathbb{P})$ and the square domain $G = [0, L]^2 \subset \mathbb{R}^2$. Let $u : G \times \Omega \to \mathbb{R}$ denote a random field over $(\Omega, \mathcal{F}, \mathbb{P})$ satisfying Assumptions* (A1) *and* (A2) *such that for $\mathbb{P}$-almost all $\omega \in \Omega$, the function $u(\cdot, \omega) : G \to \mathbb{R}$ is twice continuously differentiable. For each $\omega \in \Omega$, denote the nodal domains of $u(\cdot, \omega)$ by $N^\pm(\omega) \subset G$ and denote their cubical approximations (as in Definition 1.1) by $Q_M^\pm(\omega)$.*

*For $x = (x_1, x_2)$ and $\delta > 0$ such that $J = [x_1, x_1 + \delta] \times [x_2, x_2 + \delta] \subset G$, consider the following events in $\mathcal{F}$:*

(i) *let $E_B(x, \delta)$ denote the set of all $\omega \in \Omega$ for which $u(\cdot, \omega)$ exhibits at least one of the seven sign patterns in Figure 6 at the nine dyadic points $d_{1,\ell,m}$ for $\ell, m = 0, 1, 2$;*

(ii) *let $E_I(x, \delta)$ denote the set of all $\omega \in \Omega$ for which $u(\cdot, \omega)$ exhibits at least one of the seventeen sign patterns in Figures 8 and 10 at the nine dyadic points $d_{1,\ell,m}$ for $\ell, m = 0, 1, 2$.*

*Assume that there exist positive constants $\mathcal{C}_1$ and $\mathcal{C}_2$ such that*

$$\text{(13)} \qquad \mathbb{P}(E_B(x,\delta)) \leq \mathcal{C}_1 \cdot \delta^3 \quad \text{and} \quad \mathbb{P}(E_I(x,\delta)) \leq \mathcal{C}_2 \cdot \delta^4$$

*for all $x \in G$ and $\delta > 0$ for which $J$ lies in $G$.*

*Then for every discretization size $M$, the probability that the homologies of $N^{\pm}(\omega)$ and $Q_M^{\pm}(\omega)$ coincide satisfies*

$$\text{(14)} \qquad \mathbb{P}\{H_*(N^{\pm}) \cong H_*(Q_M^{\pm})\} \geq 1 - \frac{24\mathcal{C}_1 L^3 + 20\mathcal{C}_2 L^4}{3M^2}.$$

PROOF. Fix a discretization size $M \in \mathbb{N}$. Then the discretization points subdivide the square $G$ into $M^2$ subsquares, each with side length $\delta = L/M$. Let $J$ denote one of these subsquares which lies at the boundary. In order to estimate the probability of $J$ not being $B$-admissible, we must consider dyadic subdivisions, as in Definition 3.1. At the $n$th subdivision level, there are $4^n$ dyadic subsquares with sides of length $\delta/2^n$. According to (13), the probability that a particular one of these dyadic subsquares contains one of the forbidden sign configurations of Figure 6 is bounded by $\mathcal{C}_1 \cdot (\delta/2^n)^3$. Summing over all dyadic subsquares and subdivision levels $n \in \mathbb{N}_0$, we then obtain

$$\mathbb{P}\{J \text{ is not } B\text{-admissible}\} \leq \sum_{n=0}^{\infty} \sum_{k=0}^{4^n-1} \mathcal{C}_1 \cdot \left(\frac{\delta}{2^n}\right)^3 = \frac{2\mathcal{C}_1 L^3}{M^3}.$$

Now, assume that $J$ is contained in the interior of $G$. As above, (13) implies that

$$\text{(15)} \qquad \mathbb{P}\{J \text{ is not } I\text{-admissible}\} \leq \sum_{n=0}^{\infty} \sum_{k=0}^{4^n-1} 5\mathcal{C}_2 \cdot \left(\frac{\delta}{2^n}\right)^4 = \frac{20\mathcal{C}_2 L^4}{3M^4}.$$

Note the additional factor of 5 due to the fact that for each of the dyadic subsquares, we must also allow for the excluded sign configurations in the four shifted squares; see Definition 3.6.

According to Proposition 3.7, if $\omega \in \Omega$ is such that the homologies of $N^{\pm}(\omega)$ and $Q_M^{\pm}(\omega)$ differ, then one of the following statements must be true:

- the function $u(\cdot,\omega)$ vanishes at one of the gridpoints or at one of the dyadic subpoints of the grid subsquares;
- the function $u(\cdot,\omega)$ has a double zero;
- one of the subsquares $J$ which lies at the boundary of $G$ is not $B$-admissible;
- one of the interior subsquares $J$ is not $I$-admissible.

Due to Assumptions (A1) and (A2), the first two events have probability zero and we therefore finally obtain

$$1 - \mathbb{P}\{H_*(N^\pm) \cong H_*(Q_M^\pm)\} \leq (M^2 - (M-2)^2) \cdot \frac{2\mathcal{C}_1 L^3}{M^3} + (M-2)^2 \cdot \frac{20\mathcal{C}_2 L^4}{3M^4}$$

$$\leq \frac{24\mathcal{C}_1 L^3 + 20\mathcal{C}_2 L^4}{3M^2},$$

which completes the proof of the theorem. □

The above result could be simplified for the situation of doubly-periodic functions. If the domain of such a function is identified with a torus and we are interested in the homology of the nodal domains as subsets of this torus, then we can reformulate the theorem by employing $I$-admissibility alone. Furthermore, due to the periodicity of the domain, we can replace the additional factor of 5 in (15) by 3 since it now suffices to check each square, as well as its down-shift and its right-shift. Thus for doubly periodic functions, the final probability bound reads

$$\mathbb{P}\{H_*(N^\pm) \cong H_*(Q_M^\pm)\} \geq 1 - \frac{4\mathcal{C}_2 L^4}{M^2},$$

where $\mathcal{C}_2$ is as in (13). Finally, as will be demonstrated in the next subsection, the bounds in (13) are optimal for doubly periodic random functions.

3.4. *Random periodic functions.* As in the one-dimensional case, we close this section by providing explicit results for the important special case of doubly periodic functions, that is, for classical random Fourier series.

ASSUMPTION 3.9. We consider random Fourier series on $G = [0, L]^2$ of the form

$$
\begin{aligned}
u(x,\omega) = \sum_{k,\ell=0}^{\infty} a_{k,\ell} \cdot \Bigg( &g_{k,\ell,1}(\omega) \cos \frac{2\pi k x_1}{L} \cos \frac{2\pi \ell x_2}{L} \\
&+ g_{k,\ell,2}(\omega) \cos \frac{2\pi k x_1}{L} \sin \frac{2\pi \ell x_2}{L} \\
&+ g_{k,\ell,3}(\omega) \sin \frac{2\pi k x_1}{L} \cos \frac{2\pi \ell x_2}{L} \\
&+ g_{k,\ell,4}(\omega) \sin \frac{2\pi k x_1}{L} \sin \frac{2\pi \ell x_2}{L} \Bigg)
\end{aligned}
$$
(16)

which satisfy the following:

(a) The random variables $g_{k,\ell,m}$ in (16) are defined over a probability space $(\Omega, \mathcal{F}, \mathbb{P})$ and are independent and normally distributed with mean 0 and variance 1.

(b) there exist positive $k_1, \ell_1 \in \mathbb{N}$ and nonnegative $k_2, \ell_2 \in \mathbb{N}_0$ which satisfy $k_1 \neq k_2$ and $\ell_1 \neq \ell_2$, as well as $k_1^2 + \ell_1^2 \neq k_2^2 + \ell_2^2$, such that both $a_{k_1,\ell_1}$ and $a_{k_2,\ell_2}$ are nonzero. In addition, suppose that

$$\sum_{k,\ell=0}^{\infty} (k^6 + \ell^6) a_{k,\ell}^2 < \infty.$$

As before, the summability condition in Assumption 3.9(b) is related to smoothness properties of the random field $u$. In fact, it guarantees that the function $u(\cdot, \omega)$ has continuous partial derivatives up to order two for $\mathbb{P}$-almost all $\omega \in \Omega$. Moreover, if we define

(17) $$A_{p,q} = \sum_{k,\ell=0}^{\infty} k^{2p} \ell^{2q} a_{k,\ell}^2,$$

then we can readily show that

$$\mathbb{E}\|D_{x_1}^p D_{x_2}^q u\|_{L^2(G)}^2 = (2\pi)^{2p+2q} \cdot L^{2-2p-2q} \cdot A_{p,q},$$

where $\mathbb{E}$ denotes the expected value of a random variable over $(\Omega, \mathcal{F}, \mathbb{P})$.

As in the one-dimensional case, doubly periodic functions simplify our discussion, in that their spatial covariance function is shift-invariant. More precisely, the spatial covariance function of (16) which is given by

(18) $$R(x,y) = \mathbb{E} u(x) u(y) = \sum_{k,\ell=0}^{\infty} a_{k,\ell}^2 \cdot \cos \frac{2\pi k (x_1 - y_1)}{L} \cdot \cos \frac{2\pi \ell (x_2 - y_2)}{L}.$$

Adopting this setting, we obtain the following result.

THEOREM 3.10 (Random doubly periodic functions in 2D). *Consider the random Fourier series $u$ defined in (16) and suppose that Assumption 3.9 is satisfied. Let $M$ denote an arbitrary positive integer and let $Q_M^{\pm}(\omega)$ denote the cubical approximations of the random nodal domains $N^{\pm}(\omega)$ of $u(\cdot, \omega)$ (as in Definition 1.1) which are used for the computation of the homologies of $N^{\pm}(\omega)$. Then the probability that the homology of the random nodal domains $N^{\pm}(\omega)$ is computed correctly with the discretization of size $M$ satisfies*

(19) $$\mathbb{P}\{H_*(N^{\pm}) \cong H_*(Q_M^{\pm})\} \geq 1 - \frac{1067 \pi^2}{18 M^2} \cdot \frac{(A_{2,0} + A_{1,1} + A_{0,2})^2}{A_{0,0}^{1/2} A_{0,1}^{1/2} A_{1,0}^{1/2} A_{1,1}^{1/2}}$$
$$+ O\left(\frac{1}{M^3}\right),$$

*where $A_{p,q}$ was defined in (17).*

TABLE 1
*Covariance matrix spectral information for the first shown sign pattern. For the second pattern, we must replace $R_{0,\ell}$ by $R_{\ell,0}$*

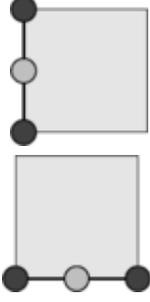

$$\det C(\delta) = \frac{R_{0,1}(R_{0,0}R_{0,2} - R_{0,1}^2)}{64} \cdot \delta^6 + O(\delta^7)$$

$$\lambda_1(\delta) = \frac{R_{0,0}R_{0,2} - R_{0,1}^2}{96R_{0,0}} \cdot \delta^4 + O(\delta^5)$$
$$\lambda_2(\delta) = \frac{R_{0,1}}{2} \cdot \delta^2 + O(\delta^3)$$
$$\lambda_3(\delta) = 3R_{0,0} + O(\delta)$$

$$v_1(\delta) = \tfrac{1}{\sqrt{6}} \cdot (1, -2, 1)^t + O(\delta)$$

TABLE 2
*Covariance matrix spectral information for the sign pattern shown*

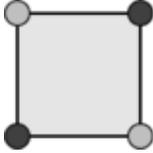

$$\det C(\delta) = R_{0,0}R_{0,1}R_{1,0}R_{1,1} \cdot \delta^8 + O(\delta^9)$$

$$\lambda_1(\delta) = \frac{R_{1,1}}{4} \cdot \delta^4 + O(\delta^5)$$
$$\lambda_2(\delta) = \frac{2R_{0,1}R_{1,0}}{R_{0,1} + R_{0,1} + |R_{0,1} - R_{0,1}|} \cdot \delta^2 + O(\delta^3)$$
$$\lambda_3(\delta) = \frac{R_{0,1} + R_{0,1} + |R_{0,1} - R_{0,1}|}{2} \cdot \delta^2 + O(\delta^3)$$
$$\lambda_4(\delta) = 4R_{0,0} + O(\delta)$$

$$v_1(\delta) = \tfrac{1}{2} \cdot (1, -1, 1, -1)^t + O(\delta)$$

PROOF. According to our assumptions, the random variable $u(x, \cdot) : \Omega \to \mathbb{R}$ is normally distributed with mean 0 and variance $\sum_{k,\ell=0}^{\infty} a_{k,\ell}^2 \neq 0$ for each $x \in G$ and this implies (A1). Furthermore, (A2) follows from [1], Theorem 3.2.1. Thus in order to apply Theorem 3.8, we need only verify (13).

As in the proof of Theorem 2.7, we employ Proposition 4.1 with $n = 3$ or $n = 4$ to estimate the probability of the various sign patterns, if one considers a square with side length $\delta$, as in (13). Since this can basically be done as in the one-dimensional case, we have only collected the necessary spectral information for the various covariance matrices in Tables 1–6. In each of these tables, the expansions of the determinant of the covariance matrix $C(\delta)$, the expansions of the eigenvalues and the expansion for the eigenvector corresponding to $\lambda_1(\delta)$ are given. The constants $R_{p,q}$ which appear in these

TABLE 3
*Covariance matrix spectral information for the sign pattern shown*

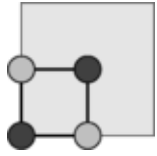

$$\det C(\delta) = \frac{R_{0,0} R_{0,1} R_{1,0} R_{1,1}}{256} \cdot \delta^8 + O(\delta^9)$$

$$\lambda_1(\delta) = \frac{R_{1,1}}{64} \cdot \delta^4 + O(\delta^5)$$
$$\lambda_2(\delta) = \frac{R_{0,1} R_{1,0}}{2(R_{0,1} + R_{0,1} + |R_{0,1} - R_{0,1}|)} \cdot \delta^2 + O(\delta^3)$$
$$\lambda_3(\delta) = \frac{R_{0,1} + R_{0,1} + |R_{0,1} - R_{0,1}|}{8} \cdot \delta^2 + O(\delta^3)$$
$$\lambda_4(\delta) = 4R_{0,0} + O(\delta)$$

$$v_1(\delta) = \tfrac{1}{2} \cdot (1, -1, 1, -1)^t + O(\delta)$$

TABLE 4
*Covariance matrix spectral information for the first sign pattern shown. For the second pattern, we must replace $R_{0,\ell}$ by $R_{\ell,0}$ and $R_{k,0}$ by $R_{0,k}$*

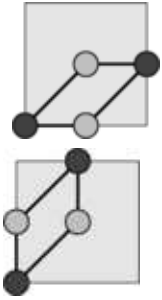

$$\det C(\delta) = \frac{R_{0,1} R_{1,0}(R_{0,0} R_{1,1} + R_{0,0} R_{2,0} - R_{1,0}^2)}{256} \cdot \delta^8 + O(\delta^9)$$

$$\lambda_1(\delta) = \frac{R_{0,0} R_{1,1} + R_{0,0} R_{2,0} - R_{1,0}^2}{64 R_{0,0}} \cdot \delta^4 + O(\delta^5)$$
$$\lambda_2(\delta) = \frac{R_{0,1} R_{1,0}}{2R_{0,1} + 4R_{1,0} + 2\sqrt{R_{0,1}^2 + 4R_{1,0}^2}} \cdot \delta^2 + O(\delta^3)$$
$$\lambda_3(\delta) = \frac{R_{0,1} + 2R_{1,0} + \sqrt{R_{0,1}^2 + 4R_{1,0}^2}}{8} \cdot \delta^2 + O(\delta^3)$$
$$\lambda_4(\delta) = 4R_{0,0} + O(\delta)$$

$$v_1(\delta) = \tfrac{1}{2} \cdot (1, -1, 1, -1)^t + O(\delta)$$

formulas arise in the expansion of the spatial covariance function (18). More precisely, if we set

$$r(d) = R(y + d, y) = \sum_{k,\ell=0}^{\infty} a_{k,\ell}^2 \cdot \cos \frac{2\pi k d_1}{L} \cdot \cos \frac{2\pi \ell d_2}{L},$$

then for $d \to 0$, the even function $r$ can be expanded as

$$r(d) = R_{0,0} - \tfrac{1}{2} \cdot (R_{1,0} d_1^2 + R_{0,1} d_2^2)$$
$$+ \tfrac{1}{24} \cdot (R_{2,0} d_1^4 + 6 R_{1,1} d_1^2 d_2^2 + R_{0,2} d_2^4) + O(|d|^6),$$

where

$$R_{k,\ell} = \left(\frac{2\pi}{L}\right)^{2k+2\ell} \cdot A_{k,\ell}.$$

In addition, we can easily verify that all of the leading coefficients in the expansions of the tables are nonzero due to Assumption 3.9(b).

By combining the information in Tables 1 and 2 with Proposition 4.1 and the definition of $E_B(x, \delta)$ in (13), we obtain the estimate

$$\mathbb{P}(E_B(x,\delta)) \leq 6 \cdot \frac{\delta^3}{128\pi} \cdot \frac{R_{0,0}R_{0,2} - R_{0,1}^2}{R_{0,0}^{3/2} R_{0,1}^{1/2}} + 6 \cdot \frac{\delta^3}{128\pi} \cdot \frac{R_{0,0}R_{2,0} - R_{1,0}^2}{R_{0,0}^{3/2} R_{1,0}^{1/2}}$$

TABLE 5
*Covariance matrix spectral information for the first sign pattern shown. For the second pattern, we must replace $R_{0,\ell}$ by $R_{\ell,0}$ and $R_{k,0}$ by $R_{0,k}$*

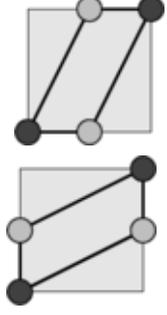

$$\det C(\delta) = \frac{R_{0,1}R_{1,0}(4R_{0,0}R_{1,1} + R_{0,0}R_{2,0} - R_{1,0}^2)}{64} \cdot \delta^8 + O(\delta^9)$$

$$\lambda_1(\delta) = \frac{4R_{0,0}R_{1,1} + R_{0,0}R_{2,0} - R_{1,0}^2}{64 R_{0,0}} \cdot \delta^4 + O(\delta^5)$$
$$\lambda_2(\delta) = \frac{R_{0,1}R_{1,0}}{2R_{0,1} + R_{1,0} + \sqrt{4R_{0,1}^2 + R_{1,0}^2}} \cdot \delta^2 + O(\delta^3)$$
$$\lambda_3(\delta) = \frac{2R_{0,1} + R_{1,0} + \sqrt{4R_{0,1}^2 + R_{1,0}^2}}{4} \cdot \delta^2 + O(\delta^3)$$
$$\lambda_4(\delta) = 4R_{0,0} + O(\delta)$$

$$v_1(\delta) = \tfrac{1}{2} \cdot (1, -1, 1, -1)^t + O(\delta)$$

TABLE 6
*Covariance matrix spectral information for the sign pattern shown*

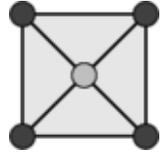

$$R^* = R_{0,0}(R_{0,2} + 2R_{1,1} + R_{2,0}) - (R_{0,1} + R_{1,0})^2$$

$$\det C(\delta) = \frac{R_{0,1}R_{1,0}R_{1,1}R^*}{64} \cdot \delta^{12} + O(\delta^{13})$$

$$\lambda_1(\delta) = \frac{R^*}{80 R_{0,0}} \cdot \delta^4 + O(\delta^5)$$
$$\lambda_2(\delta) = \frac{R_{1,1}}{4} \cdot \delta^4 + O(\delta^5)$$
$$\lambda_3(\delta) = \frac{2R_{0,1}R_{1,0}}{R_{0,1} + R_{1,0} + |R_{0,1} - R_{1,0}|} \cdot \delta^2 + O(\delta^3)$$
$$\lambda_4(\delta) = \frac{R_{0,1} + R_{1,0} + |R_{0,1} - R_{1,0}|}{2} \cdot \delta^2 + O(\delta^3)$$
$$\lambda_5(\delta) = 5R_{0,0} + O(\delta)$$

$$v_1(\delta) = \frac{1}{2\sqrt{5}} \cdot (1, 1, 1, 1, -4)^t + O(\delta)$$

$$+ 2 \cdot \frac{\delta^4}{12\pi^2} \cdot \frac{R_{1,1}^{3/2}}{R_{0,0}^{1/2} R_{0,1}^{1/2} R_{1,0}^{1/2}}$$

$$\leq \frac{3\delta^3}{64\pi} \cdot \frac{R_{0,2} R_{1,0}^{1/2} + R_{2,0} R_{0,1}^{1/2}}{R_{0,0}^{1/2} R_{0,1}^{1/2} R_{1,0}^{1/2}} + \frac{\delta^4}{6\pi^2} \cdot \frac{R_{1,1}^{3/2}}{R_{0,0}^{1/2} R_{0,1}^{1/2} R_{1,0}^{1/2}}$$

$$= \frac{3\pi^2 \delta^3}{8L^3} \cdot \frac{A_{0,2} A_{1,0}^{1/2} + A_{2,0} A_{0,1}^{1/2}}{A_{0,0}^{1/2} A_{0,1}^{1/2} A_{1,0}^{1/2}} + \frac{8\pi^2 \delta^4}{3L^4} \cdot \frac{A_{1,1}^{3/2}}{A_{0,0}^{1/2} A_{0,1}^{1/2} A_{1,0}^{1/2}}$$

$$\leq \frac{41\pi^2 \delta^3}{12L^3} \cdot \frac{(A_{2,0} + A_{1,1} + A_{0,2})^{3/2}}{A_{0,0}^{1/2} A_{0,1}^{1/2} A_{1,0}^{1/2}}.$$

Note that the three terms after the first inequality correspond to the bounds for the two patterns in Table 1 and the pattern in Table 2 (in that order); the additional factors 6, 6 and 2 correspond to the multiplicities of each pattern. For example, the first vertical pattern must be tested twice (depending on whether the interior point is positive or negative) on each of the three vertical lines.

From the above estimate, we can readily see that the first part of (13) is satisfied with

$$\mathcal{C}_1 = \frac{41\pi^2}{12L^3} \cdot \frac{(A_{2,0} + A_{1,1} + A_{0,2})^2}{A_{0,0}^{1/2} A_{0,1}^{1/2} A_{1,0}^{1/2} A_{1,1}^{1/2}}.$$

Similarly, by combining the information in Tables 3–6 with Propositions 4.1 and 4.2 (for the five-point sign pattern) and recalling the definition of $E_I(x,\delta)$ in (13), we obtain the estimate

$$\mathbb{P}(E_I(x,\delta)) \leq 8 \cdot \frac{\delta^4}{192\pi^2} \cdot \frac{R_{1,1}^{3/2}}{R_{0,0}^{1/2} R_{0,1}^{1/2} R_{1,0}^{1/2}} + 8 \cdot \frac{\delta^4}{192\pi^2} \cdot \frac{(R_{2,0} + R_{1,1})^{3/2}}{R_{0,0}^{1/2} R_{0,1}^{1/2} R_{1,0}^{1/2}}$$

$$+ 8 \cdot \frac{\delta^4}{192\pi^2} \cdot \frac{(R_{0,2} + R_{1,1})^{3/2}}{R_{0,0}^{1/2} R_{0,1}^{1/2} R_{1,0}^{1/2}} + 4 \cdot \frac{\delta^4}{384\pi^2} \cdot \frac{(R_{2,0} + 4R_{1,1})^{3/2}}{R_{0,0}^{1/2} R_{0,1}^{1/2} R_{1,0}^{1/2}}$$

$$+ 4 \cdot \frac{\delta^4}{384\pi^2} \cdot \frac{(R_{0,2} + 4R_{1,1})^{3/2}}{R_{0,0}^{1/2} R_{0,1}^{1/2} R_{1,0}^{1/2}}$$

$$+ 2 \cdot \frac{\delta^4}{1024\pi^2} \cdot \frac{(R_{0,2} + 2R_{1,1} + R_{2,0})^2}{R_{0,0}^{1/2} R_{0,1}^{1/2} R_{1,0}^{1/2} R_{1,1}^{1/2}}$$

$$\leq \frac{115\delta^4}{384\pi^2} \cdot \frac{(R_{0,2} + R_{1,1} + R_{2,0})^2}{R_{0,0}^{1/2} R_{0,1}^{1/2} R_{1,0}^{1/2} R_{1,1}^{1/2}},$$

that is, we can set

$$\mathcal{C}_2 = \frac{115\pi^2}{24L^4} \cdot \frac{(A_{2,0} + A_{1,1} + A_{0,2})^2}{A_{0,0}^{1/2} A_{0,1}^{1/2} A_{1,0}^{1/2} A_{1,1}^{1/2}}.$$

Thus, (13) is also satisfied and Theorem 3.10 follows from our abstract result, Theorem 3.8. □

It is clear from the above proof that the bound in (19) is not optimal. Its form is chosen to point out the essential ingredients, but not to establish the tightest fit. We do believe, however, that it reflects the true situation fairly accurately.

3.5. *Random bivariate trigonometric polynomials.* As in the one-dimensional case, we illustrate our main theorem for random doubly periodic functions using random trigonometric polynomials. For this, let $N \geq 3$ denote an arbitrary integer. If, in the situation of Theorem 3.10, we choose the coefficients as $a_{k,\ell} = 1$ for $1 \leq k, \ell \leq N$ and $a_{k,\ell} = 0$ otherwise, then (16) represents a bivariate random trigonometric polynomial. Now, let $A_\ell$ be defined as in (17). We can then easily show that

$$\frac{(A_{2,0} + A_{1,1} + A_{0,2})^2}{A_{0,0}^{1/2} A_{0,1}^{1/2} A_{1,0}^{1/2} A_{1,1}^{1/2}} = \frac{1}{900} \cdot (46N^2 + 51N - 7)^2 \sim \frac{529}{225} \cdot N^4.$$

This implies that in order for the homology computation to be accurate with high confidence, we must choose the discretization size $M$ in such a way that

$$M \sim N^2 \qquad \text{for } N \to \infty,$$

in contrast to the one-dimensional case. Numerically validated computational results which provide more insight into the asymptotic behavior of $M$ will be presented in [12].

**4. Probabilistic tools.** The central probabilistic tools used in the results of the previous sections are concerned with the asymptotic behavior of sign-distribution probabilities of parameter-dependent Gaussian random variables. More precisely, let $T(\delta) = (T_1(\delta), \ldots, T_n(\delta))^t \in \mathbb{R}^n$ denote a one-parameter family of $\mathbb{R}^n$-valued random Gaussian variables over a probability space $(\Omega, \mathcal{F}, \mathbb{P})$, indexed by $\delta > 0$, and choose a sign sequence $(s_1, \ldots, s_n) \in \{\pm 1\}^n$. We are then interested in the asymptotic behavior of the probability

(20) $\qquad P(\delta) = \mathbb{P}\{s_j T_j(\delta) \geq 0 \text{ for all } j = 1, \ldots, n\}.$

To this end, we derive two results. While the first result establishes the precise asymptotic behavior of $P(\delta)$ as $\delta \to 0$, the second result furnishes an upper bound under considerably weakened hypotheses.

4.1. *Precise asymptotics.* We begin with a result which establishes the precise asymptotic behavior of sign-distribution probabilities.

PROPOSITION 4.1. *Let $(s_1, \ldots, s_n) \in \{\pm 1\}^n$ denote a fixed sign sequence and consider a one-parameter family $T(\delta)$, $\delta > 0$, of $\mathbb{R}^n$-valued random Gaussian variables over a probability space $(\Omega, \mathcal{F}, \mathbb{P})$ which satisfies the following assumptions:*

(i) *For each $\delta > 0$, assume that the Gaussian random variable $T(\delta)$ has mean $0 \in \mathbb{R}^n$ and a positive-definite covariance matrix $C(\delta) \in \mathbb{R}^{n \times n}$, whose positive eigenvalues are given by $\lambda_1(\delta), \ldots, \lambda_n(\delta)$. The corresponding orthonormalized eigenvectors are denoted by $v_1(\delta), \ldots, v_n(\delta)$.*

(ii) *There exists a vector $\bar{v}_1 = (\bar{v}_{11}, \ldots, \bar{v}_{1n})^t \in \mathbb{R}^n$ such that $v_1(\delta) \to \bar{v}_1$ as $\delta \to 0$ and $s_j \cdot \bar{v}_{1j} > 0$ for all $j = 1, \ldots, n$.*

(iii) *The quotient $\lambda_1(\delta)/\lambda_k(\delta)$ converges to 0 as $\delta \to 0$ for all $k = 2, \ldots, n$.*

*Then the probability $P(\delta)$ defined in (20) satisfies*

$$(21) \qquad \lim_{\delta \to 0} P(\delta) \cdot \sqrt{\frac{\det C(\delta)}{\lambda_1(\delta)^n}} = \frac{\Gamma(n/2)}{2 \cdot \pi^{n/2} \cdot (n-1)!} \cdot \left| \prod_{j=1}^n \bar{v}_{1j} \right|^{-1}.$$

*Recall that $\Gamma(1/2) = \pi^{1/2}$, $\Gamma(1) = 1$ and $\Gamma(t+1) = t\Gamma(t)$ for $t > 0$.*

PROOF. Define the diagonal matrix $S = (s_i \delta_{ij})_{i,j=1,\ldots,n}$, where $\delta_{ij}$ denotes the Kronecker delta, and let $Z_+ = \{z \in \mathbb{R}^n : z_j \geq 0 \text{ for } j = 1, \ldots, n\}$. Finally, let

$$D(\delta) = \lambda_1(\delta) \cdot SC(\delta)^{-1}S.$$

Using the density of the Gaussian distribution of $T(\delta)$ according to [3], Theorem 30.4, which exists since $C(\delta)$ is positive-definite, in combination with a simple rescaling, the probability in (20) can be rewritten as

$$P(\delta) = \frac{(2\pi)^{-n/2}}{\sqrt{\det C(\delta)}} \cdot \int_{Z_+} e^{-z^t SC(\delta)^{-1} Sz/2} \, dz = \sqrt{\frac{\lambda_1(\delta)^n}{\pi^n \det C(\delta)}} \cdot \int_{Z_+} e^{-z^t D(\delta) z} \, dz.$$

According to our assumptions, the eigenvalues $\mu_1(\delta), \ldots, \mu_n(\delta)$ of the matrix $D(\delta)$ are given by

$$\mu_1(\delta) = 1 \quad \text{and} \quad \mu_k(\delta) = \frac{\lambda_1(\delta)}{\lambda_k(\delta)} \qquad \text{for } k = 2, \ldots, n,$$

with corresponding orthonormalized eigenvectors $w_k(\delta) = Sv_k(\delta)$ for $k = 1, \ldots, n$. Now, let $B(\delta)$ denote the orthogonal matrix with columns $w_1(\delta), \ldots, w_n(\delta)$ and introduce the change of variables $z = B(\delta)\zeta$. Moreover, let

$$Z(\zeta_1, \delta) = \left\{ (\zeta_2, \ldots, \zeta_n) : \sum_{k=1}^n \zeta_k w_k(\delta) \in Z_+ \right\} \subset \mathbb{R}^{n-1}$$

and define

$$I(\zeta_1, \delta) = \int_{Z(\zeta_1,\delta)} \exp\left(-\sum_{k=2}^{n} \mu_k(\delta)\zeta_k^2\right) d(\zeta_2, \ldots, \zeta_n).$$

Due to (ii) and the definition of the signs $s_k$, the eigenvector $w_1(\delta)$ has strictly positive components for all sufficiently small $\delta > 0$ and therefore the identity

$$z^t D(\delta) z = \sum_{k=1}^{n} \mu_k \zeta_k^2$$

implies that

(22)
$$\int_{Z_+} e^{-z^t D(\delta)z}\, dz = \int_{B(\delta)^{-1} Z_+} \exp\left(-\sum_{k=1}^{n} \mu_k \zeta_k^2\right) d\zeta$$
$$= \int_0^\infty I(\zeta_1, \delta) e^{-\zeta_1^2}\, d\zeta_1.$$

From the definition of $I(\zeta_1, \delta)$, we can easily deduce that

$$I(\zeta_1, \delta) = \zeta_1^{n-1} \cdot \int_{Z(1,\delta)} \exp\left(-\zeta_1^2 \cdot \sum_{k=2}^{n} \mu_k(\delta)\eta_k^2\right) d(\eta_2, \ldots, \eta_n)$$

and this representation implies that

(23) $\quad I(\zeta_1, \delta) \le \zeta_1^{n-1} \cdot \mathrm{vol}_{n-1}(Z(1,\delta)) \qquad$ for all $\zeta_1 > 0$ and $\delta > 0$.

Again, according to (ii), the $(n-1)$-dimensional volume of the simplex $Z(1,\delta)$ converges to the $(n-1)$-dimensional volume of the simplex

$$\widetilde{Z} = \{z \in Z_+ : (z - S\bar{v}_1, S\bar{v}_1) = 0\} \subset \mathbb{R}^n,$$

which can be computed as

$$\mathrm{vol}_{n-1}(\widetilde{Z}) = \frac{1}{(n-1)!} \cdot \left|\prod_{j=1}^{n} \bar{v}_{1j}\right|^{-1}.$$

Now, let $\zeta_1 > 0$ be arbitrary, but fixed. Note that since we did not make any assumptions about the asymptotic behavior of the eigenvectors $w_2(\delta), \ldots, w_n(\delta)$ for $\delta \to 0$, the sets $Z(1,\delta)$ do not have to converge. Yet, (ii) yields the existence of a compact subset $K \subset \mathbb{R}^{n-1}$ such that $Z(1,\delta) \subset K$ for all sufficiently small $\delta > 0$. Due to (iii), the integrand in (22) converges to 1 uniformly on $K$ and we therefore have

$$\lim_{\delta \to 0} I(\zeta_1, \delta) = \zeta_1^{n-1} \cdot \mathrm{vol}_{n-1}(\widetilde{Z}) \qquad \text{for all } \zeta_1 > 0.$$

Due to (23) and $\text{vol}_{n-1}(Z(1,\delta)) \to \text{vol}_{n-1}(\widetilde{Z})$, we can now apply the dominated convergence theorem to pass to the limit $\delta \to 0$ in (22), which leads to

$$\lim_{\delta \to 0} \int_{Z_+} e^{-z^t D(\delta) z} \, dz = \text{vol}_{n-1}(\widetilde{Z}) \cdot \int_0^\infty \zeta_1^{n-1} e^{-\zeta_1^2} \, d\zeta_1 = \text{vol}_{n-1}(\widetilde{Z}) \cdot \frac{\Gamma(n/2)}{2}.$$

This completes the proof of the proposition. $\square$

The above result is the main probabilistic tool for our one-dimensional results in Section 2. In this situation, Proposition 4.1 is used to study random variables in $\mathbb{R}^3$. We note that for three-dimensional random variables, it is, in fact, possible to explicitly compute the probability $P(\delta)$ for each $\delta > 0$ as a function of the entries in the covariance matrix $C(\delta)$. This has been demonstrated in [15] using the inversion formula ([5], Formula (29.3)) and could be used to study the asymptotic behavior of $P(\delta)$ as $\delta \to 0$. However, for the results of Section 3, we must also consider higher-dimensional random vectors and the method in [15] no longer applies.

4.2. *Upper bounds.* While the above result establishes the precise asymptotic behavior of the sign-distribution probability $P(\delta)$, the necessary assumptions for Proposition 4.1 are not satisfied in all instances discussed in Section 3. Yet, in these cases, it suffices to obtain an upper bound on the probability $P(\delta)$. Such an upper bound is the subject of the following result, which can be derived under fairly weak assumptions.

PROPOSITION 4.2. *Let $(s_1, \ldots, s_n) \in \{\pm 1\}^n$ denote a fixed sign sequence and consider a one-parameter family $T(\delta)$, $\delta > 0$, of $\mathbb{R}^n$-valued random Gaussian variables over a probability space $(\Omega, \mathcal{F}, \mathbb{P})$ which satisfies the following assumptions:*

(i) *For each $\delta > 0$, the Gaussian random variable $T(\delta)$ has vanishing mean $0 \in \mathbb{R}^n$ and a positive-definite covariance matrix $C(\delta) \in \mathbb{R}^{n \times n}$ whose positive eigenvalues are given by $\lambda_1(\delta), \ldots, \lambda_n(\delta)$. Corresponding orthonormalized eigenvectors are denoted by $v_1(\delta), \ldots, v_n(\delta)$.*

(ii) *There exists a vector $\bar{v}_1 = (\bar{v}_{11}, \ldots, \bar{v}_{1n})^t \in \mathbb{R}^n$ such that $v_1(\delta) \to \bar{v}_1$ as $\delta \to 0$ and $s_j \cdot \bar{v}_{1j} > 0$ for all $j = 1, \ldots, n$.*

*Then the probability $P(\delta)$ defined in (20) satisfies*

$$(24) \qquad \limsup_{\delta \to 0} P(\delta) \cdot \sqrt{\frac{\det C(\delta)}{\lambda_1(\delta)^n}} \leq \frac{\Gamma(n/2)}{2 \cdot \pi^{n/2} \cdot (n-1)!} \cdot \left| \prod_{j=1}^n \bar{v}_{1j} \right|^{-1}.$$

PROBABILISTIC VALIDATION OF HOMOLOGY COMPUTATIONS 39PROOF. We adopt the notation used in the proof of Proposition 4.1. We then have

$$P(\delta) = \frac{(2\pi)^{-n/2}}{\sqrt{\det C(\delta)}} \cdot \int_{Z_+} e^{-z^t SC(\delta)^{-1} Sz/2} \, dz$$
$$= \sqrt{\frac{\lambda_1(\delta)^n}{\pi^n \det C(\delta)}} \cdot \int_0^\infty I(\zeta_1, \delta) e^{-\zeta_1^2} \, d\zeta_1,$$

where

$$I(\zeta_1, \delta) = \zeta_1^{n-1} \cdot \int_{Z(1,\delta)} \exp\left(-\zeta_1^2 \cdot \sum_{k=2}^n \mu_k(\delta) \eta_k^2\right) d(\eta_2, \ldots, \eta_n)$$
$$\leq \zeta_1^{n-1} \cdot \mathrm{vol}_{n-1}(Z(1,\delta)).$$

It was shown in the proof of Proposition 4.1 that $\mathrm{vol}_{n-1}(Z(1,\delta)) \to \mathrm{vol}_{n-1}(\widetilde{Z})$ is a consequence of (ii), and this immediately implies (24). $\square$

Note that Proposition 4.2 does not require any assumption of the eigenvalue $\lambda_1(\delta)$. In particular, $\lambda_1(\delta)$ need not be the smallest positive eigenvalue of $C(\delta)$, in contrast to Proposition 4.1. We also note that we could easily modify the above result to yield lower bounds for the asymptotic behavior of $P(\delta)$. Since these are not needed for our applications, we refrain from formulating them here.

**Acknowledgment.** The authors would like to thank the anonymous referee for his careful reading of the manuscript and useful suggestions.## REFERENCES

[1] ADLER, R. J. (1981). *The Geometry of Random Fields.* Wiley, Chichester. MR0611857
[2] ADLER, R. J. and TAYLOR, J. E. (2007). *Random Fields and Geometry.* Springer, New York.
[3] BAUER, H. (1996). *Probability Theory.* de Gruyter, Berlin. MR1385460
[4] BHARUCHA-REID, A. T. and SAMBANDHAM, M. (1986). *Random Polynomials.* Academic Press, Orlando. MR0856019
[5] BILLINGSLEY, P. (1995). *Probability and Measure*, 3rd. ed. Wiley, New York. MR1324786
[6] BLÖMKER, D., MAIER-PAAPE, S. and WANNER, T. (2001). Spinodal decomposition for the Cahn–Hilliard–Cook equation. *Comm. Math. Phys.* **223** 553–582. MR1866167
[7] BLÖMKER, D., MAIER-PAAPE, S. and WANNER, T. (2007). Second phase spinodal decomposition for the Cahn–Hilliard–Cook equation. *Trans. Amer. Math. Soc.* To appear.
[8] CAHN, J. W. (1968). Spinodal decomposition. *Trans. Metallurgical Society of AIME* **242** 166–180.PROBABILISTIC VALIDATION OF HOMOLOGY COMPUTATIONS 39PROOF. We adopt the notation used in the proof of Proposition 4.1. We then have

$$P(\delta) = \frac{(2\pi)^{-n/2}}{\sqrt{\det C(\delta)}} \cdot \int_{Z_+} e^{-z^t SC(\delta)^{-1} Sz/2} \, dz$$
$$= \sqrt{\frac{\lambda_1(\delta)^n}{\pi^n \det C(\delta)}} \cdot \int_0^\infty I(\zeta_1, \delta) e^{-\zeta_1^2} \, d\zeta_1,$$

where

$$I(\zeta_1, \delta) = \zeta_1^{n-1} \cdot \int_{Z(1,\delta)} \exp\left(-\zeta_1^2 \cdot \sum_{k=2}^n \mu_k(\delta) \eta_k^2\right) d(\eta_2, \ldots, \eta_n)$$
$$\leq \zeta_1^{n-1} \cdot \mathrm{vol}_{n-1}(Z(1,\delta)).$$

It was shown in the proof of Proposition 4.1 that $\mathrm{vol}_{n-1}(Z(1,\delta)) \to \mathrm{vol}_{n-1}(\widetilde{Z})$ is a consequence of (ii), and this immediately implies (24). $\square$

Note that Proposition 4.2 does not require any assumption of the eigenvalue $\lambda_1(\delta)$. In particular, $\lambda_1(\delta)$ need not be the smallest positive eigenvalue of $C(\delta)$, in contrast to Proposition 4.1. We also note that we could easily modify the above result to yield lower bounds for the asymptotic behavior of $P(\delta)$. Since these are not needed for our applications, we refrain from formulating them here.

**Acknowledgment.** The authors would like to thank the anonymous referee for his careful reading of the manuscript and useful suggestions.

## REFERENCES

[1] ADLER, R. J. (1981). *The Geometry of Random Fields.* Wiley, Chichester. MR0611857
[2] ADLER, R. J. and TAYLOR, J. E. (2007). *Random Fields and Geometry.* Springer, New York.
[3] BAUER, H. (1996). *Probability Theory.* de Gruyter, Berlin. MR1385460
[4] BHARUCHA-REID, A. T. and SAMBANDHAM, M. (1986). *Random Polynomials.* Academic Press, Orlando. MR0856019
[5] BILLINGSLEY, P. (1995). *Probability and Measure*, 3rd. ed. Wiley, New York. MR1324786
[6] BLÖMKER, D., MAIER-PAAPE, S. and WANNER, T. (2001). Spinodal decomposition for the Cahn–Hilliard–Cook equation. *Comm. Math. Phys.* **223** 553–582. MR1866167
[7] BLÖMKER, D., MAIER-PAAPE, S. and WANNER, T. (2007). Second phase spinodal decomposition for the Cahn–Hilliard–Cook equation. *Trans. Amer. Math. Soc.* To appear.
[8] CAHN, J. W. (1968). Spinodal decomposition. *Trans. Metallurgical Society of AIME* **242** 166–180.

[9] Cahn, J. W. and Hilliard, J. E. (1958). Free energy of a nonuniform system I. Interfacial free energy. *J. Chem. Phys.* **28** 258–267.

[10] Carlsson, G. and de Silva, V. (2003). Topological approximation by small simplicial complexes. Preprint.

[11] Day, S., Kalies, W. D., Mischaikow, K. and Wanner, T. (2007). Probabilistic and numerical validation of homology computations for nodal domains. Submitted for publication.

[12] Day, S., Kalies, W. D. and Wanner, T. (2007). Homology computations of nodal domains: Accuracy estimates and validation. In preparation.

[13] de Silva, V. and Carlsson, G. (2004). Topological estimation using witness complexes. In *Eurographics Symposium on Point-Based Graphics* (M. Alexa and S. Rusinkiewicz, eds.) The Eurographics Association.

[14] Donoho, D. L. and Grimes, C. (2003). Hessian eigenmaps: Locally linear embedding techniques for high-dimensional data. *Proceedings of the National Academy of Science* **100** 5591–5596. MR1981019

[15] Dunnage, J. E. A. (1966). The number of real zeros of a random trigonometric polynomial. *Proc. London Math. Soc.* **16** 53–84. MR0192532

[16] Edelman, A. and Kostlan, E. (1995). How many zeros of a random polynomial are real? *Bull. Amer. Math. Soc.* **32** 1–37. MR1290398

[17] Farahmand, K. (1998). *Topics in Random Polynomials*. Longman, Harlow. MR1679392

[18] Gameiro, M., Mischaikow, K. and Kalies, W. (2004). Topological characterization of spatial-temporal chaos. *Phys. Rev. E* **70** 035203. MR2129999

[19] Gameiro, M., Mischaikow, K. and Wanner, T. (2005). Evolution of pattern complexity in the Cahn–Hilliard theory of phase separation. *Acta Materialia* **53** 693–704.

[20] Hyde, J. M., Miller, M. K., Hetherington, M. G., Cerezo, A., Smith, G. D. W. and Elliott, C. M. (1995). Spinodal decomposition in Fe-Cr alloys: Experimental study at the atomic level and comparison with computer models—II. Development of domain size and composition amplitude. *Acta Metallurgica et Materialia* **43** 3403–3413.

[21] Hyde, J. M., Miller, M. K., Hetherington, M. G., Cerezo, A., Smith, G. D. W. and Elliott C. M. (1995). Spinodal decomposition in Fe-Cr alloys: Experimental study at the atomic level and comparison with computer models—III. Development of morphology. *Acta Metallurgica et Materialia* **43** 3415–3426.

[22] Kahane, J.-P. (1985). *Some Random Series of Functions*, 2nd ed. Cambridge Univ. Press. MR0833073

[23] Krishan, K., Gameiro, M., Mischaikow, K. and Schatz, M. F. (2005). Homological characterization of spiral defect chaos in Rayleigh–Benard convection. Preprint.

[24] Miller, M. K., Hyde, J. M., Hetherington, M. G., Cerezo, A., Smith, G. D. W. and Elliott, C. M. (1995). Spinodal decomposition in Fe-Cr alloys: Experimental study at the atomic level and comparison with computer models—I. Introduction and methodology. *Acta Metallurgica et Materialia* **43** 3385–3401.

[25] Niyogi, P., Smale, S. and Weinberger, S. (2007). Finding the homology of submanifolds with high confidence from random samples. *Discrete and Computational Geometry*. To appear.

[26] ROWEIS, S. T. and SAUL, L. K. (2000). Nonlinear dimensionality reduction by locally linear embedding. *Science* **290** 2323–2326.
[27] SIDOROV, N., LOGINOV, B., SINITSYN, A. and FALALEEV, M. (2002). *Lyapunov–Schmidt Methods in Nonlinear Analysis and Applications*. Kluwer, Dordrecht. MR1959647
[28] TENENBAUM, J. B., DE SILVA, V. and LANGFORD, J. C. (2000). A global geometric framework for nonlinear dimensionality reduction. *Science* **290** 2319–2323.
[29] VAĬNBERG, M. M. and TRENOGIN, V. A. (1974). *Theory of Branching of Solutions of Non-Linear Equations*. Noordhoff International Publishing, Leyden.
[30] WANNER, T. (2004). Maximum norms of random sums and transient pattern formation. *Trans. Amer. Math. Soc.* **356** 2251–2279. MR2048517

DEPARTMENT OF MATHEMATICS
RUTGERS UNIVERSITY
PISCATAWAY, NEW JERSEY 08854
USA
E-MAIL: mischaik@math.rutgers.edu

DEPARTMENT OF MATHEMATICAL SCIENCES
GEORGE MASON UNIVERSITY
FAIRFAX, VIRGINIA 22030
USA
E-MAIL: wanner@math.gmu.edu